\documentclass[12pt,oneside]{amsart}
\usepackage{pslatex}
\usepackage{latexsym}
\usepackage{xy}

\xyoption{all} \CompileMatrices

\usepackage{epsfig}

\newtheorem{thm}{Theorem}[section]
\newtheorem{prop}[thm]{Proposition}
\newtheorem{lem}[thm]{Lemma}
\newtheorem{cor}[thm]{Corollary}
\newtheorem{example}[thm]{Example}

\newtheorem{dfn}[thm]{Definition}
\newtheorem{rmk}[thm]{Remark}

\begin{document}

\title{Metric properties of Outer Space}

\author{Stefano Francaviglia}
\email{s.francaviglia@sns.it}
\address{Dipartimento di Matematica Applicata ``U. Dini'' (University of Pisa)
Via Buonarroti 1/c, I-56127 Pisa (Italy) }
\author{Armando Martino}
\email{Armando.Martino@upc.edu}
\address{EPSC, edifici C1, despatx 106;
Avda. del Canal Ol\'impic, S/N 08860 Castelldefels (Spain)}

\begin{abstract}
We define metrics on Culler-Vogtmann space, which are an analogue of
the Teichm\"{u}ller metric and are constructed using stretching
factors. In fact the metrics we study are related, one being a symmetrised version of the other. We investigate the basic properties of
these metrics, showing the advantages and pathologies of both choices.

We show how to compute stretching factors between marked metric
graphs in an easy way and we discuss the behaviour of stretching
factors under iterations of automorphisms.

We study metric properties of folding paths, showing that
they are geodesic for the non-symmetric metric and, if they do not
enter the thin part of Outer space, quasi-geodesic for
the symmetric metric.
\end{abstract}

\maketitle

\tableofcontents

\section{Introduction}

Culler-Vogtmann space, or Outer Space as it is sometimes
called, has been the subject of intense study. Much of the
direction of this work has been to develop a theory for Outer Space,
and the Outer Automorphism group of a free group in an analogous way
to the theory of Teichm\"{u}ller space and the mapping class group of
a surface.

Our contribution to this effort is the study of a metric which is a
clear analogue of the Teichm\"{u}ller metric, with the goal that the
important features of both Outer Space, and the automorphisms of a
free group are captured by the geometry of this metric.

After recalling the basic definitions in section~\ref{prel}, we spend
some in section~\ref{calculating} time defining and understanding the
``one-sided'' metric, from which our metric is obtained by a
``symmetrisation''. In fact, a special case of this one-sided metric
(where the objects are a rose, and its image under an automorphism) is
a quantity that has appeared in the work of Kapovich,
\cite{MR2197815}, where it is shown that the value is computable in
double exponential time. As part of our efforts to understand our
metric, and simplify many of the proofs of its properties, we show
that the calculation is considerably simpler,
Proposition~\ref{glasses}, so that the calculation for a rose is
actually linear.

We then study the metric itself in section~\ref{metrics}, and show
that the metric topology is the same as the usual length function
topology, as well as showing that the metric is proper; closed balls
are compact in this space. This is one advantage the symmetric version
of the metric has over the unsymmetric version, since for the
one-sided metric not only are Cauchy sequences not always convergent,
but also points which {\em should} be at infinite distance, namely
points on the boundary of outer space, are actually at finite distance
from points in the interior of outer space.

Section~\ref{s5} is concerned with the connection between the geometry
of outer space and the properties of the automorphisms of a free
group. Specifically, we study the behaviour of ``folding paths'' and
their metric properties. It is fairly straightforward to show that
these paths are geodesics for the one-sided metric, but it seems to be
much more difficult to show that they are even quasi-geodesics for the
actual metric. However, these folding paths are shown to have good
properties, such as the ``4 point property'', defined in
Theorem~\ref{t4pp}.

In section~\ref{s66} we show with an example that outer space,
equipped with the symmetric metric, is not a geodesic space.
We want to stress here that such example was suggested to the authors by
Bert Wiest and Thierry Coulbois when a previous version of this paper
was posted on the arxiv.

In section~\ref{qg}, we show  that if folding paths
remain within the ``thick part'' of Outer Space, then they will be
quasi-geodesics which is a result, definitions aside, that it very
intuitive. We finish, in section~\ref{iterate} by showing that for an
automorphism of exponential growth, the map from $\mathbb{Z}$ to outer
space which sends an integer, $n$ to the $n^{th}$ iterate of a given
point under the automorphism is a quasi-isometry. Interestingly, while
this result is clearly false for automorphisms of polynomial growth,
we show that for a particular example of polynomial growth
automorphism, the folding path between the rose and a image of the
rose under an (arbitrary) iterate of the automorphism is a
quasi-geodesic with uniform constants.

\vskip\baselineskip
\paragraph*{\textsc{Acknowledgements.}}
We warmly, sincerely, and enthusiastically  thank M. Bestvina,
K. Bromberg, T. Coulbois, V. Guirardel, A. Hilion,
P. Hubert, I. Kapovich, G. Levitt, J. Los, M. Lustig, G. Th\'eret and
B. Wiest.
This research was possible thanks to the financial support of the European
Research Council (MEIF-CT-2005-010975 and MERG-CT-2007-046557.)
The
second named authors gratefully acknowledges partial support from the
MEC (Spain) through project BFM2003-06613.

\section{Preliminaries}
\label{prel}

We refer the reader to \cite{MR1950871} for an excellent survey and reference article to Culler-Vogtmann space.

Our basic objects are finite marked metric graphs of some given rank
$n$. A graph of this type is represented as a metric graph, $A$ -
that is, with a positive length assigned to every edge - and a
marking $\tau_A$ which is a homotopy equivalence from the rose with
$n$ petals, $R_n$ to A,
$$
\tau_A: R_n \to A.
$$

We shall make the standard assumption that vertices have valence at
least three. Nonetheless, we notice that it is sometimes convenient
to allow vertices of valence two. When it is clear from the contest,
we will not specify whether we use bi-valent vertices.

Two marked graphs $A$ and $B$  are equivalent if there is a
homothety, $h:A \to B$, such that the following diagram commutes up
to free homotopy,

$$
\xymatrix{
 A \ar[rr]^h &  &  B  \\
 & R_n \ar[ul]^{\tau_A} \ar[ur]_{\tau_B} \\
 }
$$

Alternatively, we could only consider metric graphs of volume 1 and
then the equivalence would be given by isometries in place of
homotheties. In either case, the resulting space of equivalence
classes is called Culler-Vogtmann space of rank $n$, or $CV_n$ (when
bi-valent vertices are allowed, two marked graphs are also equivalent
if they have a common finite subdivision.)

\begin{rmk}
In the following, if there are no ambiguities we will not distinguish
between
 a marked metric graph and its class.

When we will need to be precise we will refer to a metric graph as an
element of the unprojectivised $CV_n$, and to its class as an element
of $CV_n$.
\end{rmk}

Given any marked graph $A$,  we can look at the universal cover $T_A$
which is an $\mathbb{R}$-tree on which $\pi_1(R_n)$ acts by
isometries, via the marking $\tau_A$. (From now on, we identify the
free group of rank $n$, $F_n$, with the fundamental group of $R_n$.)
Conversely, given any minimal free action of $F_n$ by isometries on a
simplicial $\mathbb{R}$-tree, we can look at the quotient object,
which will be a graph, $A$, and produce a homotopy equivalence
$\tau_A: R_n \to A$ via the action. Equivalence of graphs in $CV_n$
corresponds to actions which are equivalent up to equivariant
homothety.

Thus, points in $CV_n$, can be thought of as equivalence classes
minimal free isometric actions on simplicial $\mathbb{R}$-trees.
Given an element $w$ of $F_n$ and a point $A$ of the unprojectivised
$CV_n$, with universal cover $T_A$ whose metric we denote by $d_A$,
we may consider,
$$
l_A(w):=\inf_{p \in T_A} d_{A}(p, wp).
$$

It is well known that this infimum is always obtained and that, for a
free action, it is non-zero for the non-identity elements of the
group. In this context, $l_A(w)$ is called the translation length of
the element $w$ in the corresponding tree and clearly depends only on
the conjugacy class of $w$ in $F_n$. Thus for any point, $A$, in
$CV_n$ we can associate the sequence $(l_A(w))_{w \in F_n}$ and it is
clear that equivalent marked metric graphs will produce two
sequences, one of which is a multiple of the other by a positive real
number (the homothety constant.) Moreover, it is also the case that
inequivalent points in $CV_n$ will produce sequences which are not
multiples of each other \cite{MR907233}. Thus, we have an
embedding of $CV_n$ into $\mathbb{R}^{F_n}/ \sim$, where $\sim$ is
the equivalence relation of homothety. The space $CV_n$ is given the
subspace topology induced by this embedding.

Finally it is clear we can realise any automorphism, $\phi$, of $F_n$
as a homotopy equivalence, also called $\phi$, of $R_n$. Thus the
automorphism group of $F_n$ acts on $CV_n$ by changing the marking.
That is, given a point $(A, \tau_A)$ of $CV_n$ the image of this
point under $\phi$ is $(A,  \tau_A \phi)$.

$$
\xymatrix{
R_n \ar[r]^{\phi} \ar@/^2pc/[rr]|{\tau_A \phi} & R_n \ar[r]^{\tau_A} & A.  \\
 }
$$

Since two automorphisms which differ by an inner automorphism always
send equivalent points in $CV_n$ to equivalent points, we actually
have an action of $Out(F_n)$ on $CV_n$, and this space is often
called {\em Outer Space} for this reason.

\section{Calculating stretching factors}
\label{calculating}

Given two marked metric graphs, $A$ and $B$ with fundamental group
free of rank $n$, we would like to compute the distance between them
and, as a first step, the ``right hand distance'' between them,
defined as follows.
\begin{dfn}[Right hand factor]
\label{rightdist} For any pair $A,B$ of marked graphs we set
$$
\Lambda_R(A,B):=\sup_{1 \neq w \in F_n} \frac{l_B(w)}{l_A(w)}.
$$
\end{dfn}

Recall that $l_A(w)$ is the translation length of the element
corresponding to $w$ in the tree $T_A$ (and hence is dependent only
on the conjugacy class of $w$). However, it is readily seen that this
translation length is the same as the length of the shortest
representative in the free homotopy class of loops in $A$ defined by
the (conjugacy class of) $w$. We note that this second definition
means that $l_A(w)$ is easy to compute given a particular $w$: we
look at the image of $w$ in $A$ via the marking and we ``cyclically
reduce'' the loop in the graph by performing free cyclic reductions
which may, of course, change the basepoint. The length of any
cyclically reduced element in this sense, calculated simply by
summing the lengths of the edges crossed, will be $l_A(w)$. We shall
also use $l_A$  to refer to the lengths of (free homotopy classes of)
loops in $A$ in the obvious way. We also note that saying a loop in
$A$ is cyclically reduced is equivalent to saying that, if we
consider the loop as a map from the circle to the graph it is an
immersion. In the same spirit, a path is reduced if it is an
immersion when considered as a map from a closed interval.

While finding lengths of elements with respect to a marked metric
graph structure is straightforward, that does not indicate how to
calculate the supremum given above. In order to do that, we need to
relate one structure to the other. One way to do this is to find an
equivariant map from $A$ to $B$, which we can simply think of as a
homotopy equivalence between the graphs, which respects the markings.
That is, a map $f$ for which the following triangle commutes up to
free homotopy,

$$
\xymatrix{
 A \ar[rr]^f &  &  B  \\
 & R_n \ar[ul]^{\tau_A} \ar[ur]_{\tau_B} \\
 }
$$

In other words, $f$ is a map homotopic to $\tau_A^{-1}$ followed by
$\tau_B$, $f \simeq \tau_B \tau_A^{-1} $. It is important to note
that this is not a graph map in that edges are not necessarily sent
to edges nor vertices to vertices. We will therefore restrict to a
particular class of maps that are more easy to handle.

\begin{dfn}[PL maps]
 A map $f:A\to B$ is a $PL$-map if it is linear on edges. More
 precisely, for each edge $e$ of $A$, if we parameterise $f|_e$ with
 the segment $[0,l_A(e)]$, then $f|_e$ has constant speed. We denote
 by $S_{f,e}$ the speed of $f|_e$ (the stretching factor of $e$.)
\end{dfn}

The stretching factor of a PL-map $f$, defined as the maximal speed
of $f$, is in fact the Lipschitz constant of $f$. We denote that
quantity by $S_f$ (the notation $L_f$ for the Lipschitz constant is
more natural but also more confusing since we already have lengths
denoted by the letter $l$)

$$S_f=\max\{S_{f,e}: e\textrm{ edge of }A\}=Lip(f)$$

In general, given $f$, there is a unique PL-map $\bar{f}$ which is
homotopic to $f$ and agrees with $f$ on vertices. It is readily
checked that
\begin{equation}
  \label{lip}
  S_{\bar f}=Lip(\bar{f})\leq Lip(f).
\end{equation}

A useful observation one can make here is that $Lip(f)$ serves as an
upper bound for $\Lambda_R(A,B)$. This is because, starting with a
loop $\gamma$ in $A$, it is clear that
$$
l_B(f(\gamma)) \leq Lip(f) l_A(\gamma).
$$
Since we can consider all loops which are cyclically reduced in $A$
this means that,
$$
l_B(w) \leq Lip(f) l_A(w),  \mbox{\rm for all } w \in F_n,$$  and we
hence proved
\begin{lem}\label{l_lip}
  For any Lipschitz map $f:A\to B$ in the homotopy class of
  $\tau_B\tau_A^{-1}$
$$\Lambda_R(A,B)=\sup_{1 \neq w \in F_n} \frac{l_B(w)}{l_A(w)} \leq Lip(f).$$
\end{lem}

 Since $f$ is arbitrary, and because of~$(\ref{lip})$,
 we can deduce that
\begin{equation}
\label{stretch} \Lambda_R(A,B)=\sup_{1 \neq w \in F_n}
\frac{l_B(w)}{l_A(w)} \leq \inf\{ S_f: f\textrm{ is PL and } f\simeq
\tau_B \tau_A^{-1}\}.
\end{equation}

It is fairly clear that the infimum on the right hand side of
equation~\ref{stretch} will be realised by an actual map.
\begin{lem}\label{lstretch}
\label{min} Let $A,B$ two marked metric graphs. Then there exists an
$f_{\infty} \simeq \tau_B \tau_A^{-1} $ such that
$$
S_{f_{\infty}}=\inf\{ S_f: f\textrm{ PL and } f \simeq \tau_B
\tau_A^{-1}\}=\inf\{Lip(f):  f\simeq\tau_B\tau_A^{-1}\}.
$$
\end{lem}
\proof For any $c$, the set of $c$-Lipschitz maps from $A$ to $B$ is
precompact by Ascoli-Arzel\'a theorem because $B$ is compact.
Therefore a sequence of maps $f_n$, whose stretching factors tend to
the infimum has a convergent sub-sequence whose limit is $f_\infty$,
and it is easily checked that $S_{f_\infty}=\inf\{S_{f_n}\}$.\qed

\begin{rmk}
  Previous lemma holds in a more general setting of spaces of length
  functions (e.g. actions on real trees.)
\end{rmk}

\begin{rmk}
Equations~\ref{lip} and~\ref{stretch},
 and Lemma~\ref{lstretch} tell us that from now
on we can, as we do, assume that any map is a PL.
\end{rmk}

Now note that there are two obstructions to making
equation~\ref{stretch} an equality. While we may realise the infimum
by a concrete map, $f$, we may still have that for a given loop
$\gamma$, not all edges of $\gamma$ may be stretched by the same
amount $S_f$. Thus we need the collection of edges which are
stretched maximally to be large enough as to contain a loop.
Furthermore, even if we have such a loop $\gamma$, the image
$f(\gamma)$ may not be cyclically reduced in $B$. However, if we have
a cyclically reduced loop, $\gamma$, in $A$, all of whose edges are
stretched by $S_f$ and such that $f(\gamma)$ is cyclically reduced in
$B$, then $\Lambda_R(A,B)=S_f$. It will turn out that there always
exists a map $f$ and a loop $\gamma$ with these properties. Before
going into details, we need some preliminaries.

\begin{dfn}
 Let $A,B$  be marked metric graphs of rank $n$.
Given a PL-map $f \simeq \tau_B \tau_A^{-1} $, we denote by
$A_{max}(f)$ the subgraph of $A$ whose edges are stretched maximally,
by $S_f$.
\end{dfn}

\begin{dfn}[Optimal maps]
  A PL-map $f \simeq \tau_B \tau_A^{-1} $ is NOT optimal if there is
  some vertex of $A_{max}(f)$ such that all edges of $A_{max}(f)$
  terminating at that vertex have $f$-image with a common terminal
  partial edge.

Otherwise $f$ is optimal.
\end{dfn}

\begin{rmk}
Using the terminology of legal and illegal turns, a PL-map is optimal
if each vertex of $A_{max}$ has at least one legal turn.
\end{rmk}

Suppose that a map $f \simeq \tau_B \tau_A^{-1}$ is not optimal. Let
$v$ be a vertex of $A_{max}(f)$ such that all edges of $A_{max}(f)$
terminating at $v$ have $f$-image with a common terminal partial
edge, say $\alpha$. Let star$(v)$ denote the set of edge emanating
from $v$. We set $N=\textrm{star}(v)\cap A_{max}$ and
$K=\textrm{star}(v)\setminus N$.

Now, let $f_t$ be the homotopy that moves $v$ backward along
$\alpha$. More precisely, we let $F:A\times[0,T]\to B$ be the
homotopy such that $f_t=F(\cdot,t):A\to B$ is the PL-map that agrees
with $f$ outside star$(v)$ and such that $f_t(v)\in\alpha$ with
$d(f_t(v),\alpha)=t$. Such a homotopy exists for small $t$. Moreover,
for small $t$ we have:
\begin{enumerate}
\item For any $e_0\in N$ and any $e_1\in K$,  $S_{f,e_1}<S_{f,e_0}$.
\item
\begin{itemize}
\item Either $S_{f_t}=S_f$ and $A_{max}(f_t)\subset A_{max}(f)$ (but
  not equal.)
\item Or $S_{f_t}<S_f$ and $A_{max}(f_t)=A_{max}(f)$.
\end{itemize}
\end{enumerate}

\begin{dfn}
Let $t_0$ be the supremum of times $t$ such that $f_t$ exists and has
the above properties. We define $Next_v(f)$ as $f_{t_0}$.
\end{dfn}

Note that $Next_v(f)$ can be defined only for non-optimal maps. We
can now prove that the inequality~\ref{stretch} is an equality, as
was first proved by Tad White.

\begin{prop}\label{p_optimal}
\label{equal}  Let $A,B$ be marked metric graphs of rank $n$. Then
there exists an $f \simeq \tau_B \tau_A^{-1}$ and a cyclically
reduced loop $\gamma$ contained in $A_{max}$, the subgraph of
maximally stretched edges of $A$, whose $f$-image is also cyclically
reduced. In particular, $\Lambda_R(A,B)=S_f$ for this map $f$.
\end{prop}

\proof By Lemma~\ref{min}, we may choose a map $f$ whose stretching
factor is minimal. Moreover, we may choose such a map with the least
number of edges in $A_{max}(f)$. Hence, $Next_v(f)$ cannot exist, and
therefore $f$ is optimal. This means that any path, $p$, in
$A_{max}(f)$ which is mapped to a reduced path by $f$ can be
continued to a longer path, which is also mapped to something
reduced. This is because the obstruction to continuing $p$ is exactly
non-optimality of $f$. Starting with a single edge, and since there
are only finitely many oriented edges in $A_{max}(f)$, we can find a
reduced path of the form $eqe$ which is mapped to a reduced path by
$f$. It is then clear that $\gamma=eq$ is a cyclically reduced loop,
which is mapped to something cyclically reduced. Moreover,
$l_B(\gamma)=S_f l_A(\gamma)$, and hence $\Lambda_R(A,B)=S_f$ as
required. \qed

\medskip

Actually, one can do better.
\begin{dfn}
  Let $f:A\to B$ be a PL-map.
  For any sub-graph $A_0$ of $A$, we define $\partial_fA_0$ the
  $f$-boundary of $A_0$ as the set of vertices $v$ of $A_0$ such that all
  edges of $A_0$ terminating at $v$ have $f$-image with a common
  terminal partial edge.
\end{dfn}

So, for example, a map is optimal if and only if
$\partial_fA_{max}=\emptyset$.

\begin{prop}
 Let $A,B$ be marked metric graphs of rank $n$.
Then there exists an $f \simeq \tau_B \tau_A^{-1}$ such that, if
 $\lambda_1>\cdots>\lambda_k$ are the stretching factors of edges, if
$A_i$ denotes the sub-graph of edges stretched by $\lambda_i$, then
for all $i$
$$\partial_f A_i\subset A_{i-1}.$$
(So, heuristically, $A_i$ is a cycle relative to $A_{i-1}$.)
\end{prop}
\proof Once one founds optimal maps as in Proposition~\ref{equal},
choose between them one that has the smallest $\lambda_2$ and $A_2$,
argue as in Proposition~\ref{equal}, and conclude inductively on
$i$.\qed

\medskip

 We note
that implicit in the proof of Proposition~\ref{equal} is a proof that
$\Lambda_R(A,B)$ is computable. Namely, the path $\gamma$ produced at
the end of the proof can be chosen minimally, and so we may assume
that it passes through each oriented edge at most once. There are
only finitely many such paths, and we may compute their lengths in
$A$ and $B$ (without reference to $f$) as well as the maximum of the
ratio of these lengths. By the Proposition, this maximum will be
exactly $\Lambda_R(A,B)$. However, the number of such $\gamma$ will
be exponential in the number of edges. We will now show that it is
always possible to find a ``less complicated'' loop $\gamma$, which
will cut down the computational complexity considerably.

We will approach this problem in two steps, and the idea of this
result is that we want to reduce the complexity of $\gamma$ as a loop
in $A$. We always have in mind an optimal map $f$, and so we will
assume that $\gamma$ lies in $A_{max}$. We shall attempt to simplify
by cutting and gluing $\gamma$ to itself. Since we will only use
edges that were already in $\gamma$, we ensure that our loops are
always contained in $A_{max}$. In order for the cutting and pasting
to result in loops which still give the value for $\Lambda_R(A,B)$,
we need to make sure that the resulting image in $B$ is cyclically
reduced. Therefore we always need to keep in mind that we are working
at two levels. On the one hand we have a loop, $\gamma$, thought of
as a map from the circle to $A$ ($A_{max}$, in fact). We then compose
this map with $f$ and the resulting loop in $B$ is an immersion. For
the first step of our result, we prove the following ``Sausages
Lemma'', which says that we may take a $\gamma$ which realises
$\Lambda_R(A,B)$ and whose shape in $A$ has in Figure~\ref{sausage
picture}.

\begin{figure}[htbp]
\begin{center}
\ \includegraphics[width=3in]{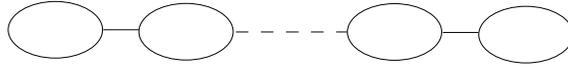}
\end{center}
\caption{Sausages} \label{sausage picture}
\end{figure}

\medskip

For any oriented path $\gamma$ we denote by $\overline{\gamma}$ its
inverse.

\begin{lem}[Sausages Lemma]
\label{sausages}  Let $A,B$ be marked metric graphs of rank $n$, and
let $f\simeq \tau_B\tau_A^{-1}$ be an optimal map. Then there exists
a loop $\gamma$ such that $l_B(\gamma)/ l_A(\gamma)=S_f =
\Lambda_R(A,B) $. In particular, $\gamma$ is cyclically reduced in
$A$ and in $B$ via $f$. Furthermore, $\gamma$ is a sausage, i.e.
$\gamma=\gamma_1 \overline{\gamma_2}$ where each $\gamma_i$ is a path
in $A$ that can be parameterised with $[0,1]$ in such a way that
\begin{itemize}
\item $\gamma_1$  and $\gamma_2$ are embeddings;
\item there exists a finite family of disjoint closed intervals $I_j\subset
  (0,1)$, each one possibly consisting of a single point, such that
  $\gamma_1(t)=\gamma_2(s)$ if and only if $t=s$
  and $t$ belongs to $\{0,1\}\cup_j I_j$.
\end{itemize}
\end{lem}
\proof The content of the result is that $\gamma=
\gamma_1\overline{\gamma_2}$ with the specified properties, since
everything else follows from Proposition~\ref{equal}. This will
follow from two sublemmas. First we establish some notation. We shall
think of $\gamma$ as a map from $S^1$ to $A$ and also, via $f$, as a
map from $S^1$ to $B$. We shall subdivide $S^1$ to give it a graph
structure and so that edges map to edges in $B$. For simplicity,
although it isn't really necessary, we shall assume that all the
vertices of $A$ map to vertices of $B$, which we can arrange after a
suitable subdivision.

\begin{figure}[ht]
\begin{center}
\ \includegraphics[width=1.5in]{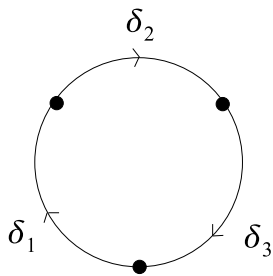}
\end{center}
\caption{Triple Points} \label{deltas}
\end{figure}

Our first sublemma says that if three distinct points in $S^1$ have
the same image in $A$, then we can choose $\gamma$ to be shorter (in
both $A$ and $B$.) To do this, we look at three points in $S^1$
mapped to the same point in $A$. Thus we decompose $\gamma$ as
$\delta_1 \delta_2 \delta_3$ as in the picture above, where the
endpoints of each $\delta_i$ map to the same point in $A$. Our first
attempt is to try to replace $\gamma$ with one of the $\delta_i$,
each of which is clearly a shorter path in $A$, and each of which
maps to a reduced path in $B$. The only way that this can fail is if
each $\delta_i$ maps to a reduced but not cyclically reduced path in
$B$. This means that we can write,
$$
\begin{array}{rcl}
\delta_1 & = & e_1 \ldots \overline{e_1} \\
\delta_2 & = & e_2 \ldots \overline{e_2} \\
\delta_3 & = & e_3 \ldots \overline{e_3}, \\
\end{array}
$$
where we are writing each $\delta_i$ as a concatenation of edges {\em
labelled} by the image of that edge in $B$. Thus we are saying that
the image of $\delta_1$ in $B$ begins with an edge $e_1$ and ends
with the inverse edge $ \overline{e_1}$. However, we know that
$\gamma$ is immersed in $B$, so that $e_1 \neq  e_2$. In particular,
this implies that the loop $\delta_1 \delta_2$ is immersed in $B$,
and we are done.

\begin{figure}[ht]
\begin{center}
\ \includegraphics[width=1.8in]{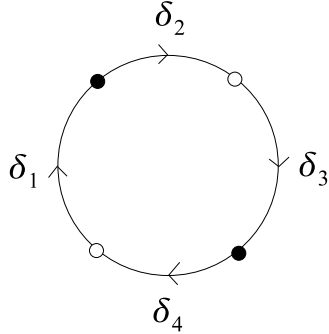}
\end{center}
\caption{Crossing Points} \label{moredeltas}
\end{figure}

For the second sublemma, we will show that we can avoid crossing
double points in $\gamma$. That is, if we can write $\gamma$ as a
concatenation $\delta_1 \delta_2 \delta_3 \delta_4$ in $S^1$ such
that the initial points of $\delta_1$ and $\delta_3$ have the same
image in $A$, and the initial points of $\delta_2$ and $\delta_4$
have the same image in $A$, then we may replace $\gamma$ by a shorter
path (shorter in both $A$ and $B$).

Now we try to replace $\gamma$ by one of the paths $\delta_i
\delta_{i+1}$ (subscripts taken modulo 4). If any of these map to
cyclically reduced loops in $B$, we are done. Otherwise, we get that,
$$
\begin{array}{rcl}
\delta_1 \delta_2 & = & e_1 \ldots \overline{e_1} \\
\delta_2 \delta_3 & = & e_2 \ldots \overline{e_2} \\
\delta_3 \delta_4 & = & e_3 \ldots \overline{e_3} \\
\delta_4 \delta_1 & = & e_4 \ldots \overline{e_4}, \\
\end{array}
$$
where, as before, this is a concatenation of edges in $S^1$ labelled
by the images in $B$. This implies that
$$\delta_i = e_i \ldots \overline{e_{i+3}},$$
with subscripts taken mod 4. Since we know that $\gamma$ is immersed
in $B$, we must have that $e_1 \neq e_3$ and $e_2 \neq e_4$. Thus it
is clear that the loop $\delta_1 \overline{\delta_3}$ is immersed in
$B$, and hence we have proven the second sublemma.

For our third and final sublemma, we wish to remove all
``bad triangles''.  
This may be slightly confusing terminology, but we wish to avoid the
situation where $\gamma$ is the concatenation of 6 paths, where alternating
paths in this decomposition are closed (and the other 3 form a, not
necessarily embedded, triangle).
Formally, let us assume that we can write
$$
\gamma= \delta_1 \delta_2 \delta_3 \delta_4 \delta_5 \delta_6, $$
where $\delta_1, \delta_3, \delta_5$ are closed paths, and show that this means
we can shorten $\gamma$. Note that if any of the paths $\delta_1, \delta_3, \delta_5$
are immersed in $B$ then we are done, simply by replacing $\gamma$. So let us assume
that none of these subpaths are immersed. Using similar arguments as before, this implies that
 $\delta_2 \delta_4 \delta_6$ is a closed path which is immersed in $B$ and we are done.

Armed with these sublemmas, we may remove all triple points,
all crossing points and all bad triangles since there are only finitely many loops less
than a given length in $A$ (or $B$). We subdivide $\gamma$ into edges and vertices,
labelled by their image in $B$. Clearly, the labelling need not be unique since $\gamma$ need
not be embedded, however since we have removed all triple points each label may occur at most
twice. If each label occurs once, $\gamma$ is embedded and we are done. Otherwise, choose an
"innermost" pair of vertices in $\gamma$ with the same label. That is, choose such a pair $u, v$ and
one of the paths between them, $\delta$ so that $\delta$ embeds in $B$ except at the endpoints.

Since we have removed all bad triangles, there are at most two innermost such pairs (in fact there are exactly two,
if we also keep track of the path between them and remember that we are assuming that $\gamma$ is not embedded).
For each innermost pair, choose a point between them (ie. on the specified path). So we now have two points on $\gamma$ and therefore
two subpaths, $\gamma_1, \gamma_2$ between them and $\gamma=\gamma_1 \overline{\gamma_2}$.
Since we have no bad triangles, both $\gamma_1$ and $\gamma_2$ must be embedded in $B$.
We have also divided $\gamma$, and hence its subpaths, according to the image in $B$ and use this
parameterisation to finish the Lemma. Namely, the disjoint intervals $I_j$ correspond to edges or vertices of $B$
which have more than one pre-image in $\gamma$. Since we have eliminated all crossing points in $\gamma$, the intervals $I_j$
appear in the same order in both $\gamma_1$ and $\gamma_2$ are we are done.
\qed

\medskip

The final step in simplifying our loop $\gamma$ is to move from a
collection of sausages to at most two.

\begin{prop}
\label{glasses}
 Let $A,B \in CV_n$,
and let $f\simeq \tau_B\tau_A^{-1}$ be an optimal map. Then there
exists a loop $\gamma$ with
$l_B(\gamma)/l_A(\gamma)=S_f=\Lambda_R(A,B)$
 so that either
\begin{enumerate}
 \item[O.] $\gamma$ is a simple closed curve in $A$,
\item[\Large$\infty$.] $\gamma$ is an embedded bouquet of two circle, i.e.
$\gamma=\gamma_1 \gamma_2$, where $\gamma_i$ are simple closed curves
which do not meet each other, except at a single point, or
\item[O$-$O.] $\gamma=\gamma_1 \gamma_3 \gamma_2\overline{\gamma_3}$, where
$\gamma_1$ and $\gamma_2$ are simple closed curves which do not meet,
and $\gamma_3$ is an embedded path that touches $\gamma_1$ and
$\gamma_2$ at their initial points only.
\end{enumerate}
In particular, there exists a finite set of loops, $\Gamma$, in $A$
so that $l_B(\gamma)/l_A(\gamma)=\Lambda_R(A,B)$ for some $\gamma \in
\Gamma$ and the set $\Gamma$ can be chosen {\em independently} of
$B$.
\end{prop}

\proof We shall start by taking the loop
$\gamma=\gamma_1\overline{\gamma_2}$ supplied by
Lemma~\ref{sausages}. If the family of intervals $\{I_j\}$ is empty,
then $\gamma$ is a simple closed curve; if it consists of a single
interval $I$ then $\gamma$ is either an embedded {\Large$\infty$}- or
O$-$O-curve, depending whether $I$ is a single point or not. In these
cases we are done.

Suppose that the family $\{I_j\}$ contains at least two intervals. We
show how to reduce to the case of only one interval. Let $[a,b]$ and
$[c,d]$ be the two extremal intervals of $\{I_j\}$; namely, such that
$0<a\leq b< c\leq d<1$ and no $I_j$ in $(0,a)\cup(d,1)$. We replace
the loop $\gamma_2$ with the following
$$
\delta_2(t)=\left\{\begin{array}{lcl}
\gamma_2(t)&\qquad& t<b\\
\gamma_1(t)&\qquad& t\in[b,c]\\
\gamma_2(t)&\qquad& t>c\end{array}\right.
$$

Note that $\delta_2$ is embedded in $A$ because
$\gamma_1(t)=\gamma_2(s)$ if and only if $t=s$ (by
Lemma~\ref{sausages}.) Also, the $f$-image of $\delta_2$ in $B$ is
reduced because of the same reason and because the $f$-images of both
$\gamma_1$ and $\gamma_2$ are reduced. The new loop
$\widetilde\gamma=\gamma_1\overline{\delta_2}$ is therefore a
sausage-loop satisfying
$l_B(\widetilde\gamma)/l_A(\widetilde\gamma)=S_f=\Lambda_R(A,B)$, and
the cardinality of the $I_j$'s is now one.\qed

\medskip

Another interesting consequence of Proposition~\ref{p_optimal} is
that $\Lambda_R$ is always defined and finite. We notice that this
can also be proved directly using the immersion of paths in the space
of geodesic currents. Indeed, the space of geodesic currents is
compact, and lengths are continuous linear functionals, so the ratio
of two length functionals always has maximum and minimum realised by
some current. In particular the maximum is finite and the minimum is
non-zero, and we have additionally proved that it is realised by a {\em rational} current.

\section{Metrics}
\label{metrics}

We are now in a position to define a metric on $CV_n$ and our
starting point will be Definition~\ref{rightdist}. In fact we have
both left hand and right hand displacements (whose existence is
guaranteed by Proposition~\ref{equal} and the preceding discussion.)

\begin{dfn}[Right hand left factors] For any pair $A,B$ of marked
metric graphs of rank $n$ we set:
$$
\Lambda_R(A,B):=\sup_{1 \neq w \in F_n} \frac{l_B(w)}{l_A(w)} \qquad
\Lambda_L(A,B):=\sup_{1 \neq w \in F_n}
\frac{l_A(w)}{l_B(w)}=\Lambda_R(B,A).
$$
\end{dfn}

\begin{rmk}
  Since $F_n$ embeds in the space of geodesic currents as a dense
  sub-space, we could equivalently define $\Lambda_R$ and $\Lambda_L$
  taking the supremum over the space of currents.
\end{rmk}

The reason that we wish to use both $\Lambda_R$ and $\Lambda_L$ is
that they are not symmetric functions and hence if we wish to define
a genuine metric on $CV_n$ we will need to use both of them. We are
now ready to define the metric on $CV_n$.

\begin{dfn}[Distance]
\label{metric} For all $A,B \in CV_n$, we define
$$
\Lambda(A,B) := \Lambda_R(A,B) \Lambda_L(A,B).
$$
The {\em distance} between $A$ and $B$ is then given by,
$$
d(A,B)= \log \Lambda(A,B).
$$
\end{dfn}

The first remark is that if we scale the length functions $l_A$ and
$l_B$ by positive numbers, $d(A,B)$ remains unchanged. So it is
well-defined on $CV_n$ with values (a priori) in $[-\infty,\infty]$.

Proposition~\ref{equal} shows in fact that $d(A,B)$ is always finite
(which is straightforward using currents,) but we still need to show
that it is indeed a metric. We begin with an elementary observation.

\begin{rmk}
Given a positive real valued function, $f$,
$$
\sup \frac{1}{f(x)}= \frac{1}{\inf f(x)}.
$$
Moreover, $\sup \frac{1}{f(x)}$ exists if and only if $\inf f(x)$
exists and is non-zero.
\end{rmk}

This has an easy but interesting consequence for us,
\begin{lem}
\label{supinf}
$$
\Lambda(A,B)=\frac{\sup_{1 \neq w \in
F}\,\frac{l_A(w)}{l_B(w)}}{\inf_{1 \neq w \in
F}\,\frac{l_A(w)}{l_B(w)}}
$$
\end{lem}
\proof Apply the previous remark to $\frac{l_A(w)}{l_B(w)}$, noting
that $\Lambda_R(A,B)$ always exists. \qed

It is now immediate that $d$ will be a non-negative function,

\begin{cor}
For all $A, B \in CV_n$, $\Lambda(A,B) \geq 1$ and hence $d(A,B) \geq
0$.
\end{cor}

Next we need to show that $d$ is only zero when the two entries are
the same point of $CV_n$.

\begin{lem}
Given $A, B \in CV_n$, $d(A,B)=0$ if and only if $A=B$.
\end{lem}
\proof Thinking of $CV_n$ as a space of length functions, it is clear
that if the two functions, $l_A$ and $l_B$ differ by a multiplicative
constant, then $\Lambda(A,B)=1$ and so $d(A,B)=0$. Conversely, if
$d(A,B)=0$ then after rescaling (by $\Lambda_R(A,B)$) we get that
$l_A=l_B$. \qed

\begin{lem}[Triangular inequality]\label{l_trian}
For all marked metric graphs $A,B,C$ of rank $n$
$$d(A,C) \leq d(A,B)+d(B,C).$$
\end{lem}
\proof For any $1\neq g \in F_n$
$$
\begin{array}{rcl}
\Lambda_R(A,B) \Lambda_R(B,C) & = & \sup_{1 \neq w \in F_n}
\frac{l_B(w)}{l_A(w)} \sup_{1 \neq w' \in F_n}
\frac{l_C(w')}{l_B(w')} \\ &  \geq & \frac{l_B(g)}{l_A(g)}
\frac{l_C(g)}{l_B(g)} \\ & = & \frac{l_C(g)}{l_A(g)}.
\end{array}
$$
Thus $\Lambda_R(A,B) \Lambda_R(B,C) \geq \Lambda_R(A,C)$. Using the
same argument for $\Lambda_L$, we have verified the triangle
inequality for $d$. \qed

Since the function $d$  is clearly symmetric, collecting previous
lemmata we have a proof of

\begin{thm}
The function $ d(A,B)=\log \Lambda(A,B) $ defines a metric on $CV_n$.
\end{thm}

\begin{rmk}
  It is straightforward that automorphisms of the free group act by
  isometries on $CV_n$ with respect to $d$.
\end{rmk}

Armed with the metric above, we clearly need to verify that the
topology it gives is the same as the one we already have on $CV_n$.

\begin{thm}[The topology]\label{t_top}
The topology induced by $d$ on CV is the usual one.
\end{thm}
\proof First of all, recall that marked metric graphs are
characterised by their translation lengths, so elements of  $CV_n$
are characterised by the projective classes of their translation
lengths.

We  show that the two topologies have the same converging sequences,
that being enough since both topologies have countable bases.

First, we show that if $d(A_k,A)\to 0$ then $A_k\to A$ in $CV_n$. If
$d(A_k,A)\to 0$, then by Lemma~\ref{supinf} the function

$$\frac{\sup}{\inf}(l_{A_k}/l_A)$$
uniformly converges to $1$. Therefore, up possibly to rescaling,
$l_{A_k}\to l_A$ pointwise, and thus $A_k\to A$ as elements of
$CV_n$.

Conversely, if $A_k\to A$ as elements of $CV_n$, then, up possibly to
rescaling, $A_k\to A$ as marked metric graphs. Therefore, there exist
$h_k\to 1$ and $h_k$-Lipschitz functions $f_k:A_k\to A$ and $g_k:A\to
A_k$ in the homotopy classes corresponding to the markings.
Therefore, Lemma~\ref{l_lip} (and its analogous for $\Lambda_L$)
implies  $d(A_k,A)\to 0$.\qed

\begin{thm}[Completeness]
  For any $X\in CV_n$, any closed $d$-ball centred at $X$ is
  compact. Whence $(CV_n,d)$ is complete.
\end{thm}

\proof Let $\{A_i\}$ be any sequence in $CV_n$ such that
$d(X,A_i)\leq e^R$. We show that it has a convergent sub-sequence. By
hypothesis we have
$$\frac{\sup}{\inf}(l_{A_i}/l_X)<R$$
and, up to possibly scaling the metric of $A_i$, we can suppose $\inf
(l_{A_i}/l_X)=1$. Therefore $\{\sup (l_{A_i}/l_X)\}$ is a bounded
sequence, and a diagonal argument now shows that, up to possibly passing to subsequences,  $l_{A_i}$ has as pointwise limit that we denote by
$l_\infty$. Since the closure of Outer Space is the space of ``very
small actions'', \cite{MR907233}, \cite{bf}, \cite{MR1341810},
$l_\infty$ corresponds to a translation length function of a minimal
isometric action of the free
group $F_n$ on an $\mathbb R$-tree. Since the infimum of functions is
upper semicontinuous, $l_\infty$ is bounded below away from zero. We
show in Lemma~\ref{l_inf}
 that this implies that the
action given by $l_\infty$ is actually free on a simplicial tree, and
corresponds therefore to a point $A$ of $CV_n$ which, by
Theorem~\ref{t_top} is the limit of $\{A_i\}$.\qed

\begin{lem}\label{l_inf}
  Let $l$ be the translation length function of a minimal isometric
  action of the free group $F_n$ on a $\mathbb R$-tree $T$. If $\inf l
  >c>0$ then $T$ is simplicial and the action is free.
\end{lem}
\proof The fact that the action is free is obvious since $l$ is
bounded below away from zero. Now suppose, by contradiction, that the
action is not simplicial. Then, there is a point $x\in T$ and a
sequence of segments $\sigma_k$, no three of them co-linear, such
that the sequence $\{s_k\}$ of their starting points converges to
$x$. Let $\widetilde R_n$ denote the universal cover of the standard
rose $R_n$ (i.e. $\widetilde R_n$ is the Cayley graph of $F_n$) with
a marked origin $O$,
 and let $f:\widetilde R_n \to T$ be a Lipschitz, $PL$-map
which is equivariant with respect to the actions of $F_n$ on $R_n$
and $T$. let $y_k\in R_n$ such that $f(y_k)=s_k$. Let $w_k\in F_n$ be
elements such that $w_k(y_k)$ stay at distance less than one from
$O$. After passing to a subsequence, we may assume that $w_k(y_k)$ is
convergent in $\widetilde R_n$, and hence that $w_k(s_k)$ is
convergent in $T$. Looking at distances in $T$ we see that,
$$
\begin{array}{rcl}
d(w_k(s_k), w_h(s_k)) &  \leq  & d( w_k (s_k), w_h(s_h)) + d(w_h(s_h), w_h(s_k)) \\
& = & d( w_k (s_k), w_h(s_h)) + d(s_h), s_k).
\end{array}
$$
Hence, from the remarks above, the translation length of ${w_h}^{-1}
w_k$ in $T$, tends to zero, as $h,k\to \infty$. Moreover, since the
no three of the $\sigma_k$'s are co-linear, the family $\{w_k\}$ is
infinite and hence ${w_h}^{-1} w_k$ cannot always equal the identity.
This contradicts the hypothesis that $l$ is bounded away from zero.
\qed

\bigskip
Since our metric $d$ is the corresponding of a symmetrised version of
the Thurston metric on Teichm\"uller space, it is natural to ask what
happens to the non-symmetric pieces.

\begin{dfn}
  Given $A\in CV_n$ we denote by $\bar A$ its representative which has
  total volume one.
\end{dfn}

\begin{dfn}[Right and left hand non-symmetric metric]
  For any $A,B\in CV_n$ we define
$$d_R:=\log(\Lambda_R(\bar A, \bar B))\qquad
d_L:=\log(\Lambda_L(\bar A, \bar B)).$$
\end{dfn}

Since $d_L(A,B)=d_R(B,A)$ we can restrict our study to the right hand
metric $d_R$. The elementary properties require some more work than
in the case of the symmetric metric.

First of all, note that $d$ is well-defined for marked metric graph,
and it is scale-invariant, so it descends to a metric on $CV_n$. This
property does not hold for $d_R$, however, which is why the
normalisation to volume one is crucial.

\begin{lem}
For any $A,B\in CV_n$ the right hand distance is non-negative and
vanishes only if $A=B$:
$$d_R(A,B)\geq 0 \qquad {\rm and} \qquad d_R(A,B)=0 \Leftrightarrow
A=B\in CV_n.$$
\end{lem}
\proof Let $f:\bar A\to \bar B$ be an optimal map (that exists by
Proposition~\ref{p_optimal}) then
\begin{equation}
\label{dR} 1={\rm vol}(\bar B)={\rm vol}({\rm Im}(F))\leq
\Lambda_R(\bar A, \bar B){\rm vol}(\bar A)= \Lambda_R(\bar A, \bar B)
\end{equation}
so $d_R(A,B)\geq0$. If, for any edge $e$ of $\bar A$ we denote by
$l_{\bar
  A}(e)$ its length (hence $\sum l_{\bar A}(e)=1$)
recalling that $S_{f,e}$ denotes the stretching factor of $e$, we
have

\begin{equation}\label{dR2}
{\rm vol}({\rm Im}(f))= \sum_{e \textrm{ edge of }A} S_{f,e} l_{\bar
A}(e) -C
\end{equation}

where $C$ is a non-negative quantity that measure overlappings of
$f$. Therefore, if $\Lambda_R(A,B)=1$, then the inequality of
$(\ref{dR})$ is an equality, and from $(\ref{dR2})$ we get
$S_{f,e}=1$ for all edges $e$, and $C=0$ which together imply that
$f$ is an isometry. Thus $\bar A=\bar B$ as marked graphs, and $A=B$
as elements of $CV_n$. \qed

\bigskip

Ordered triangular inequality is already proven in
Lemma~\ref{l_trian}, so we have proved

\begin{thm}
  The function $d_R(A,B)$ defines a non-symmetric metric on $CV_n$.
\end{thm}

As for the symmetric case, the topology induced by $d_R$ on $CV_n$ is
the usual one.
\begin{thm}[The Topology]
  For any sequence $\{A_k\}$ and $A\in CV_n$
$$d(A,A_k)\to 0 \Leftrightarrow d_R(A,A_k)\to 0 \Leftrightarrow
d_R(A_k,A)\to 0.$$
\end{thm}
Clearly if $d=d_R+d_L\to 0$ then both $d_R$ and $d_L$ go to zero.
Suppose that $d_R(A,A_k)\to 0$. Let $f_k:\bar A\to \bar A_k$ be an
optimal map. As in $(\ref{dR2})$ we have

$$
1={\rm vol}({\rm Im}(f_k))= \left(\sum_{e \textrm{ edge of }A}
S_{f_k,e}\, l_{\bar A}(e)\right)
 - C_k
$$
with $S_{f_k,e}\leq \Lambda_R(A,A_k)\to 1$ and $\sum l_{\bar
  A}(e)=1$. Which implies that $f_k$ converges to an isometry and
therefore $d(A,A_k)\to 0$. A similar argument works for when
$\Lambda_R(A_k,A)\to 1$.\qed

The first important difference between symmetric and non-symmetric
metrics is that the latter are not complete. Therefore, in general,
the fact that a sequence is a right hand Cauchy sequence does not
guarantee convergence in $CV_n$.

\begin{thm}[Incompleteness]\label{thm.4.19}
The space $(CV_n,d_R)$ is not complete. Namely there are sequences
$\{A_k\}$ such that $d_R(A_k,A_{k+m})\to 0$ as $k\to \infty$ which
have no accumulation point.
Moreover, for any $A\in CV_n$ and any $B\in \overline{CV_n}\setminus CV_n$
one has that $\Lambda_R(A,B)<\infty$.
\end{thm}

 \proof Let $A_0$ be $R_n$ the
standard $n$-petals rose with a uniform metric of volume one. Let
$A_k$ be the graph obtained by multiplying the metric of one petal by
a factor $1/k$ and normalised to have volume one. Then, a direct
calculation shows
$$\Lambda_R(A_k,A_{k+m})=\frac{((k+m)n-1)k}{(k+m)(kn-1)}$$
which goes to $1$ as $k\to\infty$. Thus, $\{A_k\}$ is a right hand
Cauchy sequence, but its only accumulation point is the standard rose
with $n-1$ petals which does not belong to $CV_n$ -- but it can be
viewed as an element of $\overline{CV_n}$.

In order to prove the second statement, one simply constructs a PL,
equivariant map from $A$ to $B$. This is guaranteed to be Lipschitz,
since $A$ is in $CV_n$ (for any choice of $B$.) Whence
$\Lambda_R(A,B)$ is bounded.
\qed

\begin{rmk}
  Theorem~\ref{thm.4.19} points out another ``pathology'' of the
  non-symmetric metrics. Indeed, consider a volume-one, marked metric graph $A$,
  and a sequence $B_k$ of volume-one, marked metric graphs such that
  $\Lambda_R(A,B_k)$ goes to infinity. This can be easily done
  using iterations of automorphisms (see for instance
  Section~\ref{iterate}.) Then, up to possibly passing to a
  subsequence, $B_k\to B$ a point in
  $\overline{CV_n}\setminus CV_n$. By Theorem~\ref{thm.4.19} we have
  $\Lambda_R(A,B)<\infty$ and $\Lambda_R(A,B_k)\to\infty$.
\end{rmk}

\medskip

On the other hand, right and left hand metrics are more deeply
related to folding procedures, this providing an easy description of
geodesics.

We note that one interesting consequence of the existence of the metric, is that one can use it to prove the Bounded Cancellation Lemma of \cite{cooper}.

The Bounded Cancellation Lemma, first proved by Cooper, is a key
result in the study of automorphisms of free groups. It has many
equivalent formulations, of which we state one.

\begin{thm}[Bounded Cancellation Lemma, \cite{cooper}]
Let $A, B$ be marked metric graphs of rank $n$, and consider $f:A \to
B$, a PL map such that $f \simeq \tau_B \tau_A^{-1} $. Let $|.|_A$
and $|.|_B$ denote the length functions of $A$ and $B$ respectively.
(Note that this is not quite the translation length, since we do not
cyclically reduce). Let $\alpha, \beta$ be loops in $A$, at a vertex
$v$, such that $|\alpha \beta|_A= |\alpha|_A + |\beta|_A$. Then, there
exists a constant $K$ depending only on $A$ and
$B$ (and not on $\alpha, \beta$) such that,
$$|f(\alpha \beta)|_B \geq |f(\alpha)|_B+ |f(\alpha)|_B - 2K.$$
We call $K$ a bounded cancellation constant for the map, $f$, which
clearly only depends on $f$ up to homotopy relative to
vertices.
\end{thm}

We observe that the existence of the bounded cancellation constant is
related to our (left) distance.

\begin{prop}
\label{short}
Given $A, B$ and $f$ as above, let $\lambda$ be the Lipschitz
constant for $f$. Then if $i$ is {\em not} a bounded cancellation
constant for $f$, we may find loops
$\alpha_i, \beta_i$ at a vertex $v$ of $A$ such that
\begin{enumerate}
\item $|\alpha_i \beta_i|_A= |\alpha_i|_A + |\beta_i|_A$
\item $|f(\alpha_i \beta_i)|_B < |f(\alpha_i)|_B+ |f(\beta_i)|_B -
  2(i-\lambda vol(A))$
\item $|f(\alpha_i)|_B \leq \lambda vol(A) + i$, $|f(\beta_i)|_B \leq
\lambda vol(A) + i$.
\end{enumerate}
 Moreover, we can ensure that $\alpha_i \beta_i$
is cyclically reduced in $A$.
\end{prop}
\proof By hypothesis, we may find loops $\alpha_i$, $\beta_i$ such
that  $|f(\alpha_i \beta_i)|_B < |f(\alpha_i)|_B+ |f(\beta_i)|_B -
2i$. This means that there is a terminal segment of $f(\alpha_i)$ cancels with an initial segment of $f(\beta_i)$ of length $i$ (though the cancellation may be longer). We can look at the
pre-image of this segment in $\alpha_i$ and $\beta_i$. Now, by adding
a segment of length not greater than $vol(A)$ to each of these
pre-images, we may replace $\alpha_i, \beta_i$ by paths which are
loops, (which we continue to call $\alpha_i, \beta_i$) so that
$\alpha_i \beta_i$ is cyclically reduced in $A$.

By construction, $f(\alpha_i)$ is a loop in $B$ which is the original
cancellation segment of length $i$, followed by a path which is the
image of something of length at most $vol(A)$. Since the image of
this terminal segment has length at most $\lambda vol(A)$, we know
that a terminal segment of $f(\alpha_i)$ of length at least $i -
\lambda vol(A)$ survives (and is a terminal segment of the original
cancellation segment). By a similar argument for $f(\beta_i)$, we may
deduce that a segment of length at least $i-\lambda vol(A)$ must
cancel in $f(\alpha_i \beta_i)$. Therefore, $|f(\alpha_i \beta_i)|_B
< |f(\alpha_i)|_B+ |f(\alpha_i)|_B - 2(i-\lambda vol(A))$.

Moreover, by construction, $|f(\alpha_i)|_B \leq \lambda vol(A) + i$,
$|f(\beta_i)|_B \leq \lambda vol(A) + i$ and we are done. \qed

\vskip\baselineskip

Now, consider two loops in $A$, $\alpha, \beta$,
which are based at the same vertex of $A$, such that $\alpha \beta$ is
cyclically reduced and $|\alpha \beta |_A = |\alpha|_A + |\beta|_A$,
and with the additional contidion that
$|\alpha|_A, |\beta|_A \leq 4 \lambda vol(A)\Lambda_L(A, B)$.
Let $$K_{\alpha,\beta}=
\frac{|f(\alpha)|_B + |f(\beta)|_B - |f(\alpha \beta)|_B}{2}$$
since there are only finitely many pairs, $\alpha,\beta$
with the above properties, we may find a maximum $K$ of the numbers
$K_{\alpha, \beta}$.
\begin{cor}
With the above notation, the number $K+\lambda vol(A)$ is a bounded
cancellation constant for $f$.
\end{cor}
\proof Recall that $$
\frac{1}{\Lambda_L(A,B)}=\inf_w\frac{||w||_B}{||w||_A}$$
 and that $||w||\leq |w|$ with
equality if and only if $w$ is cyclically reduced.
In particular, whenever $\alpha\beta$ is cyclically reduced, we have
$$\frac{|f(\alpha\beta)|_B}{|\alpha\beta|_A}\geq
\frac{||f(\alpha\beta)||_B}{||\alpha\beta||_A}\geq
\frac{1}{\Lambda_L(A,B)}.$$
By
Proposition~\ref{short}, if $i$ is not a bounded cancellation constant
for $f$, we may find $\alpha, \beta$ such that $\alpha\beta$ is
cycliclally reduced, the cancellation in
$f(\alpha \beta)$ is greater than $i-\lambda vol(A)$, and
$|f(\alpha)|_B,|f(\beta)|_B\leq \lambda vol(A)+i$.

So we get $|f(\alpha\beta)|_B\leq 4\lambda vol(A)$ and
$$|f(\alpha\beta)|_B\Lambda_L(A,B)\geq|\alpha\beta|_A$$
whence $|\alpha\beta|_A\leq 4\lambda vol(A)\Lambda(A,B)$ and thus
$K_{\alpha,\beta}\leq K$.

Since the cancellation in
$f(\alpha \beta)$ is greater than $i-\lambda vol(A)$
$$|f(\alpha\beta)|_B\leq|f(\alpha)|_B+|f(\beta)|_B-2(i-\lambda vol(A))$$
whence $i-\lambda vol(A)\leq K$.\qed

\section{Folding paths and geodesics}\label{s5}

In this section we study properties of geodesics and metric
properties of folding paths for the symmetric and the non-symmetric
metrics.

The following lemma provides an easy characterisation of geodesics
\begin{lem}\label{l3_triang_eq}
  Let $\gamma$ be a continuous path from an interval $[a,b]$ to a
  (possibly non-symmetric) metric space. If for any three points
  $x<y<z\in[a,b]$ $\gamma$
  realises the triangular equality
$$d(\gamma(x),\gamma(y))+d(\gamma(y),\gamma(z))=d(\gamma(x),\gamma(z)),$$
then $\gamma$ is geodesic.
\end{lem}
\proof Given a subdivision $a=t_0<t_1<\dots<t_n=b$ of $[a,b]$, the
sum $\sum_{i=1}^n d(\gamma(t_{i-1}),\gamma(t_i))$ approximates the
length of $\gamma$ as the subdivision is finer and finer. By the
triangular equality we get
$$
\sum_{i=1}^nd(\gamma(t_{i-1}),\gamma(t_i)) =
d(\gamma(t_0),\gamma(t_2))+\sum_{i=3}^nd(\gamma(t_{i-1}),\gamma(t_i))
$$
and inductively we conclude that $\gamma$ is rectifiable and that its
length realises the distance between $\gamma(a)$ and $\gamma(b)$.
\qed

\begin{cor}\label{c3_triang_eq}
  Let $A_t, t\in[a,b]$ denote a continuous path in $CV_n$. Suppose
  that for each $x,y,z\in[a,b]$ there is a loop $\gamma$ which is
  maximally stretched both from $A_x$ to from $A_y$ and $A_y$ to
  $A_z$. More precisely, suppose that
$$
\max_w\frac{l_{A_y}(w)}{l_{A_x}(w)} =
\frac{l_{A_y}(\gamma)}{l_{A_x}(\gamma)} \qquad
\max_w\frac{l_{A_z}(w)}{l_{A_y}(w)} =
\frac{l_{A_z}(\gamma)}{l_{A_y}(\gamma)}.
$$

Then $A_t$ is a $d_R$-geodesic.
\end{cor}
\proof It is immediate to check that $A_t$ realises the (oriented)
triangular equality.\qed

\medskip

The very same argument gives the following
\begin{cor}
  Let $A_t, t\in[a,b]$ denote a continuous path in $CV_n$. Suppose
  that for each $x,y,z\in[a,b]$ there are loops $\gamma$ and $\eta$ which are
  respectively maximally and minimally stretched both from $A_x$ to
  from $A_y$ and $A_y$ to
  $A_z$. More precisely, suppose that
$$
\max_w\frac{l_{A_y}(w)}{l_{A_x}(w)} =
\frac{l_{A_y}(\gamma)}{l_{A_x}(\gamma)} \qquad
\max_w\frac{l_{A_z}(w)}{l_{A_y}(w)} =
\frac{l_{A_z}(\gamma)}{l_{A_y}(\gamma)};
$$
$$
\min_w\frac{l_{A_y}(w)}{l_{A_x}(w)} =
\frac{l_{A_y}(\eta)}{l_{A_x}(\eta)} \qquad
\min_w\frac{l_{A_z}(w)}{l_{A_y}(w)} =
\frac{l_{A_z}(\eta)}{l_{A_y}(\eta)}.
$$

Then $A_t$ is a $d$-geodesic.
\end{cor}

\begin{rmk}\label{r6_5.4}
  Since $d=d_R+d_L$, a path is $d$-geodesic if and only if it is both
  $d_R$- and $d_L$-geodesic.
\end{rmk}

We are now ready to construct $d_R$-geodesics using scalings and
folding paths.

\begin{thm}[Right hand geodesics]\label{t3_drgeod}
For each $A,B$ in $CV_n$ there is a $d_R$-geodesic path between them,
that is to say a continuous path $t\mapsto A_t$ such that
$d_R(A,A_t)=t$ and $A_{d_R(A,B)}=B$.
\end{thm}
\proof Recall that $\bar A$ and $\bar B$ denote the volume-one
representatives in their respective projective classes. Let $f:\bar
A\to \bar B$ be an optimal map, let $\gamma\subset \bar A_{max}$ be a
path realising $\Lambda_R(\bar A,\bar B)$. Namely, $\gamma$ is a
geodesic in $\bar A$ (i.e. a reduced path) whose $f$-image is
geodesic (i.e. reduced) in $\bar B$, and such that $\gamma$ is
uniformly stretched by $f$ exactly by $\Lambda_R(\bar A,\bar B)$. The
existence of such $f$ and $\gamma$ is ensured by
Proposition~\ref{p_optimal}.

Let $A'$ be the marked metric graph obtained by $\bar A$ by shrinking
each edge so that it is stretched by $f$ exactly by $\Lambda_R(\bar
A,\bar B)$, and let $A_0$ the graph homothetic to $A'$ so that
$\Lambda_R(A_0,\bar B)=1$. We still denote by $f$ the induced map
$f:A_0\to \bar B$.

Note that we still have that $\gamma$ is a reduced loop in $A_0$
whose $f$-image is reduced, and that it realises the maximal
stretching factor $\Lambda_R(A_0,\bar B)=1$. Also, note that now $f$
stretches each edges of $A_0$ exactly by $1$ (that is to say, $f$ is
an isometry o edges.)

We describe now a folding procedure that will produce our geodesic.
The idea is that we never touch $\gamma$, so that it will realise the
maximally stretched loop between any two points of the folding, so
that we can invoke Corollary~\ref{c3_triang_eq}.

First, we subdivide --- allowing valence-two vertices --- both $A_0$
and $\bar B$ so that $f$ is simplicial ({\em i.e.} vertices to
vertices, edges to edges.) For each vertex $v$ of $A_0$ and $t\geq 0$
let $\sim_{t,v}$ be the equivalence relation on $A_0$ defined by:
$$x\sim_{t,v} y$$ if and only if $f(x)=f(y)$ and both $x$ and $y$
lies at distance less or equal than $t$ from $v$. Let $\sim_t$ be the
union of all relations $\sim_{t,v}$ as $v$ varies on the set all
vertices of $A_0$. For $t\geq 0$ we define
$$A_t:=A_0/\sim_{t,v}$$
we denote by $p_t$ the projection $A_0\to A_t$, and we denote by
$f_t$ the map $A_t\to \bar B$ induced by $f$, which is well-defined
since $x\sim_{t,v} y$ implies $f(x)=f(y)$ .

For small times $t$, $A_t$ is obtained from $A_0$ just identifying
germs of edges having the same image under $f$ (local folding.) Let
$t_1$ be the smallest time $t$ such that a pair of edges of $A_0$ is
completely identified in $A_t$.

Our first claim is that, for $t\in [0,t_1]$, $A_t$ is a metric graph
and that $f_t$ is an homotopy equivalence, whence $A_t$ is a marked
metric graph. The fact that $A_t$ is a graph is because for any
segment $\sigma$ in $\bar B$, $f^{-1}(\sigma)$ is a finite union of
segments, and therefore $A_t$ is the result of identifications of a
finite number of segments. The fact that $f_t$ is a homotopy
equivalence follows from the fact that $f$ factorises as
$$f:A_0\stackrel{p_t}{\to}A_t\stackrel{f_t}{\to}\bar B$$
and from the fact that $(p_t)_*:\pi_1(A_0)\to\pi_1(A_t)$ is
surjective.

Our second claim is now that $\gamma$ realises both
$\Lambda_R(A_0,A_t)$ and $\Lambda_R(A_t,\bar B)$. First note that, as
the $f$-image of $\gamma$ is geodesic, then also its $p_t$-image is.
Thus $\Lambda_R(A_0,A_t)$ is greater or equal to the ratio
$l_{A_t}(\gamma)/l_{A_0}(\gamma)$, which is one because $p_t$ is a
local isometry on edges, fact that also implies that
$\Lambda_R(A_0,A_t)\leq1$. Thus $\Lambda_R(A_0,A_t)=1$ and it is
realised by $\gamma$. A similar argument shows that $\Lambda(A_t,
\bar B)=1$ is realised by $\gamma$.

We argue now by induction. As above, we define relations
$\sim_{t-t_1,v}$ for each vertex $v$ of $A_1$, and $\sim_{t-t_1}$ as
their union. For $t\geq t_1$, we set $A_t:=A_{t_1}/\sim_{t-t_1}$, we
let $p_t:A_0\to A_t$ be the projection, and $f_t$ be map induced by
$f$.

 As above, it is easy to check that we have that
 $\Lambda_R(A_0,A_t)=\Lambda_R(A_t,\bar B)=1$ are both realised by
 $\gamma$.

 Our third claim is that such a process ends in a finite time.
Indeed, since our folding is isometric on edges, for $t<s$ we can
bound below the difference of volumes
 $${\rm vol}(A_{t})-{\rm vol}(A_{s})$$
by $s-t$.

So we must stop at a time, say $\bar t$. Since we stopped, at each
vertex the folding relations are trivial, but this simply means that
$f_{\bar t}$ is an isometry.

Summarising, we have constructed a path $A_t$ for $t\in[0,\bar t]$
with the property that, $A_0$ is in the class of $A'$ as element of
$CV_n$, $A_{\bar t}=\bar B$ is in the class of $B$ as element of
$CV_n$ and for each $t$ $\Lambda_R(A_0,A_t)$ and $\Lambda_R(A_t,\bar
B)$ are realised by the same $\gamma$. This last property does not
change if we rescale each $A_t$ to its volume-one multiple $\bar
A_t$. Therefore, for each $t$ we have
$$d_R(A',B)=d_R(A',A_t)+d_R(A_t,B).$$

Now, note that for any $0\leq s<t\leq\bar t$, if we construct a
folding path from $A_s$ to $A_t$ following the above rules, we find
exactly the restriction of the folding path we build so far.
Therefore, the path $A_t$ from $A_0$ to $B$ realises the triangular
equality, and is therefore $d_R$-geodesic by
Lemma~\ref{l3_triang_eq}.

The shrinking procedure from $A$ to $A'$ also realises the triangular
equality because everything is shrank and $\gamma$ is not touched.
Finally, if we consider a point $X$ between $A$ and $A'$ and a point
$Y$ on the geodesic between $A'$ and $B$, again we have that every
loop is stretched less than $\Lambda_R(A,B)$ and $\gamma$ is
stretched exactly by $\Lambda_R(A,B)$. In conclusion, $\gamma$ always
realises the maximum stretching factor between any two points in the
path we constructed. Such a path is then $d_R$-geodesic by
Corollary~\ref{c3_triang_eq}. \qed

\medskip
Since it is of independent interest, we formalise the precise
definition and notation the folding procedure described in the proof
of Theorem~\ref{t3_drgeod}.

\begin{dfn}[Fast folding paths and turns]
\label{fastfold} Let $A,B$ be two marked metric graphs, let $f:A\to
B$ be an optimal map, and let $A_0,\bar B$ as in the proof of
Theorem~\ref{t3_drgeod}.

A {\em fast folding path} is a path $t\mapsto A_t$ constructed
following the procedure described in the proof of
Theorem~\ref{t3_drgeod}.

A fast folding path comes with the simplicial subdivisions and the
sequence of times $0=t_0<t_1<\dots<\bar t$ such that in each
$[t_i,t_{i+1}]$ a whole segment is identified.

A {\em turn} $\tau$ at a time $t$ of a fast folding path is a pair of
edges having a common end-point and whose germs are identified for
$t'>t$. We say that the turn is folded, or that $\tau$ is a folding
turn.
\end{dfn}

\begin{rmk}
  The folding path we constructed in the proof of
  Theorem~\ref{t3_drgeod} is not unique in general, as in general we
  can start folding at many different vertices. This shows that
  $d_R$-geodesics between points of $CV_n$ are not unique.
\end{rmk}

We analyse now the local structure of geodesics in the PL structure
of $CV_n$.
\begin{dfn}[Simplices of Outer Space]
  A simplex of $CV_n$ is a sub-set of $CV_n$ consisting of all marked
  metric graphs
  with fixed topological type and marking.
\end{dfn}

Given a marked graph with edges $e_1,\dots,e_k$, the corresponding
simplex $\sigma$ is identified with the positive cone of $\mathbb
R^k$ just by assigning the metric, i.e. a length for each edge:
$$A\in\sigma \longleftrightarrow (l_A(e_1),\dots,l_A(e_k)).$$

Similarly, we can assign to each loop, its counting vector. Namely,
for a loop $\xi$ let $\xi(e_i)$ be the number of occurrences of the
edge $e_i$ in $\xi$; then
$$\xi\mapsto (\xi(e_1),\dots,\xi(e_k)).$$

This viewpoint generalises immediately to the setting of geodesic
currents (see~\cite{MR2197815}, \cite{MR2216713},
\cite{francaviglia}) and in fact it is in that setting that linear
structures arises naturally. Nevertheless, since the use of currents
is not strictly necessary for our purposes, we stick to the world of
loops.

The local linear structures of $CV_n$ and the space of loops have as
consequence that we can handle length as a linear function
$$L_A(\xi)=\langle A,\xi\rangle:=\langle(l_A(e_1),\dots,l_A(e_k)),(\xi(e_1),\dots,\xi(e_k))\rangle$$
where the last scalar product is the standard one of $\mathbb R^k$.

\begin{prop}\label{p6_5.9}
  Segments in simplices of $CV_n$ are $d_R$- and $d_L$-, whence $d$-, geodesics.
\end{prop}
\proof Let $A,B$ marked metric graphs in the same simplex. Let $\xi$
be a loop that realises $\sup_w l_B(w)/l_A(w)$. The segment between
$A$ and $B$ is parameterised by $A_t=(1-t)A+tB$ (as vectors of
$\mathbb R^k$.) For any $0\leq s<t\leq 1$ we have
$$\frac{l_{A_t(w)}}{l_{A_s}(w)}
=\frac{(1-t)\langle A, w\rangle + t\langle B, w\rangle} {(1-s)\langle
A, w\rangle + s\langle B, w\rangle}= \frac{(1-t) + t(\langle B,
w\rangle/\langle A, w\rangle)} {(1-s) + s(\langle B, w\rangle/\langle
A, w\rangle)}.
$$
The function
$$x\mapsto \frac{1-t+tx}{1-s+sx}$$
is monotone increasing for $t>s$. Thus, for any $s<t$ the stretching
factor $l_{A_t(w)}/l_{A_s}(w)$ is maximal on $\xi$. The thesis now
follows from Corollary~\ref{c3_triang_eq}. \qed

\begin{example}
  Two points of the same simplex are connected by several geodesics.
\end{example}
\proof In $CV_2$, consider the simplex of the trivalent graph with a
disconnecting edge (i.e. a O-O graph.) Let $A$ be the vector
$(1,1,1)$ where the middle coordinate is referred to the
disconnecting edge. Let $B=(\lambda,1,\lambda^{-1})$ with $\lambda
>1$, and let $c=(\lambda,1,1)$.
Let $\gamma_1$ be the segment between $A$ and $B$. Let $\gamma_2$ be
the union of the segment between $A$ and $C$ and the one between $C$
and $B$.  Using Corollary~\ref{c3_triang_eq} it is readily checked
that $\gamma_1$ and $\gamma_2$ are different geodesics between $A$
and $B$.\qed

\section{The symmetric metric is not geodesic}\label{s66}
In this section, we describe an example of two points in $CV_2$ which
are not connected by a $d$-geodesic.

This example is due to Bert Wiest and Thierry
Coulbois.

\medskip

\begin{figure}[htbp]
\begin{center}
\includegraphics[width=5in]{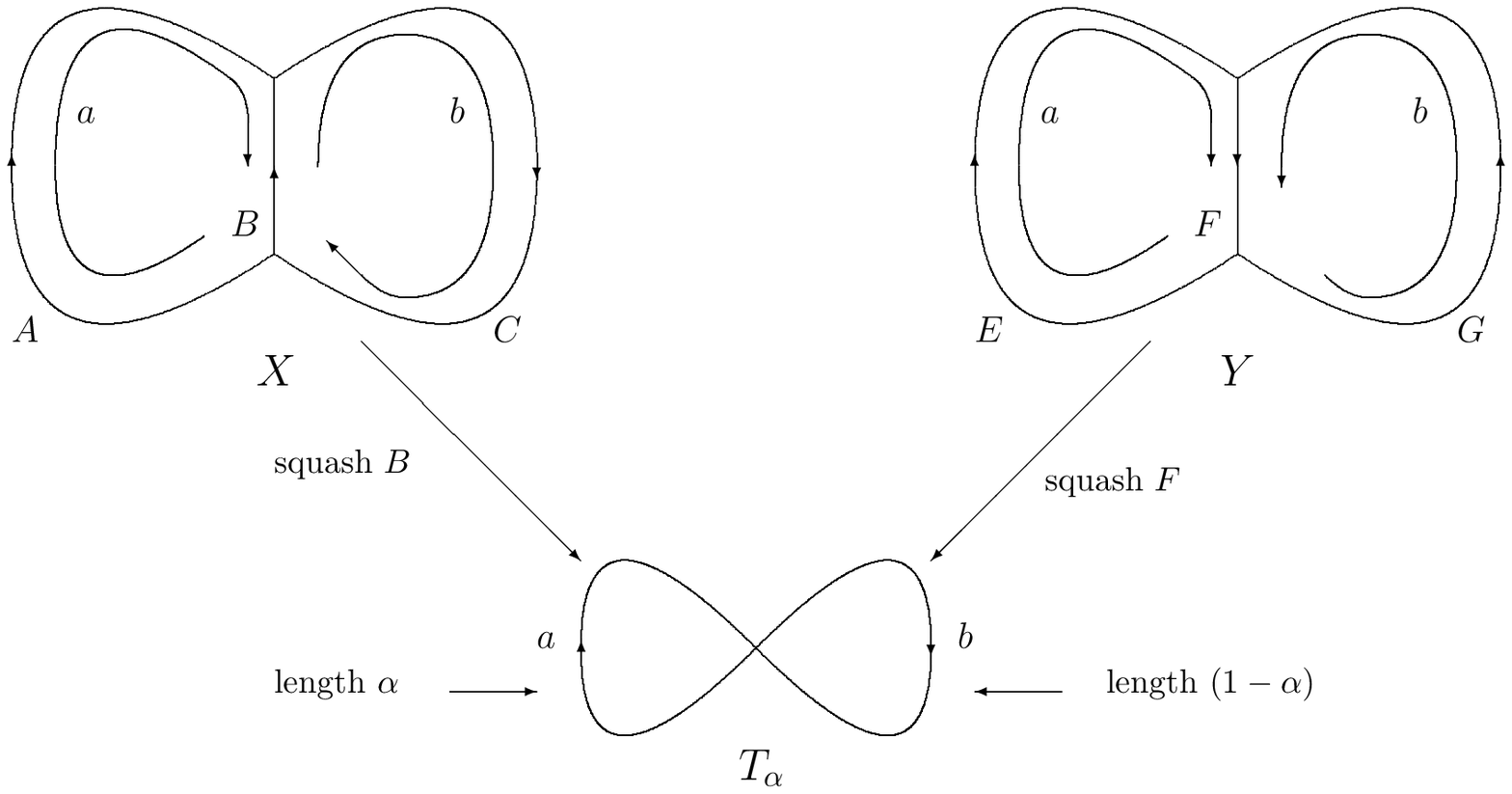}
\end{center}
\caption{The graphs $X,Y$ and $T_\alpha$} \label{BT}
\end{figure}

Consider the outer space in rank two, with graphs normalised to have
volume one, and where we denote the generators
of the free group of rank two by $a$ and $b$.
Consider two simplex of maximal dimension in $CV_2$
corresponding to graphs without disconnecting edges (theta-graphs)
such that they touch along a $1$-dimensional simplex corresponding to
a rose with two petals. Let $X$ and $Y$ be two points metric graphs,
one in each simplex, as shown in Figure~\ref{BT}.

Since each
$1$-simplex disconnects $CV_2$, any
path between $X$ and $Y$ must cross the edge common to the two
simplices. We parameterise such edge by a number $\alpha$, so that in the graph
$T_\alpha$ the petal corresponding to $a$ has length $\alpha$, and the
one corresponding to $b$ has length $1-\alpha$ (see figure~\ref{BT}.)

By proposition~\ref{p6_5.9} a $d_R$-geodesic between $X$ and $Y$
reduces to the union of two segments $XT_\alpha$ and $T_\alpha Y$,
for some $\alpha$. By Remark~\ref{r6_5.4}, if there is a $d$-geodesic
between $X$ and $Y$, there exists $\alpha$ such that $XT_\alpha\cup
T_\alpha Y$ is both $d_R$- and $d_L$-geodesic. It is readily checked
that  $XT_\alpha\cup T_\alpha Y$ is $d_R$-geodesic if and only if
there is a loop which is maximally stretched from $X$ to $T_\alpha$,
from $T_\alpha$ to $Y$ and from $X$ to $Y$ (so that the triangular
inequality become equality.) The same holds for $d_L$.

We choose now $X$ and $Y$ in a suitable way, we compute the $\alpha$ so
that $XT_\alpha\cup T_\alpha Y$ is $d_R$-geodesic and we show that
for such $\alpha$ $XT_\alpha\cup T_\alpha Y$ is not $d_L$-geodesic.

We choose $X$ and $Y$ in a symmetric way with respect the common edge:
$$X: \qquad A=1/6 \qquad B=1/3 \qquad C=1/2$$
$$Y: \qquad E=1/2 \qquad F=1/3 \qquad G=1/6$$

We compute now the right factors $\Lambda_R(X,T_\alpha)$ and
$\Lambda_R(X,Y)$. By Proposition~\ref{glasses} we have to check only
the lengths of the loops $AB, BC, AC$.

\begin{center}
\begin{tabular}[h]{|l|r|r|r|}
  \hline
Loop in $X$  & $AB$ & $BC$ & $AC$\\
\hline
Length in $X$ & $1/2$ & $5/6$ & $2/3$\\
\hline
Length in $T_\alpha$ & $\alpha$ & $1-\alpha$ & $1$\\
\hline
$l_{T_\alpha}/l_X$ & $2\alpha$ & $6(1-\alpha)/5$ & $3/2$\\
\hline
Corresponding loop in $Y$ & $EF$ & $GF$ & $EFGF$\\
\hline
Length in $Y$ & $5/6$ & $1/2$ & $4/3$\\
\hline
$l_Y/l_X$ & $5/3$ & $3/5$ & $2$\\
\hline
Loop maximally stretched from $X$ to $Y$ & ~& ~& $*$\ \ \ \\
\hline
\end{tabular}
\end{center}

It follows that $AC$ must be the maximally stretched also from $X$ to
$T_\alpha$, whence we get
$$3/2\geq2\alpha \qquad \textrm{and} \qquad 3/2\geq 6(1-\alpha)/5$$
that is $$\alpha\leq 3/4.$$

We compute now $\Lambda_R(T_\alpha,Y)$. By Proposition~\ref{glasses}
we have only to check the loops $a, b, ab, ab^{-1}$

\begin{center}
\begin{tabular}[h]{|l|r|r|r|r|}
  \hline
Loop in $T_\alpha$  & $a$ & $b$ & $ab$ & $ab^{-1}$ \\
\hline
Length in $T_\alpha$ & $\alpha$ & $1-\alpha$ & $1$ & $1$\\
\hline
Corresponding loop in $Y$ & $EF$ & $GF$ & $EG^{-1}$& $EFGF$ \\
\hline
Length in $Y$ & $5/6$ & $1/2$ & $2/3$ & $4/3$\\
\hline
$l_Y/l_{T_\alpha}$ & $5/6\alpha$ & $1/2(1-\alpha)$ & $2/3$ & $4/3$\\
\hline
\end{tabular}
\end{center}

thus, since $ab^{-1}$ must be the maximally stretched loop, we get
$$4/3\geq 5/6\alpha\qquad \textrm{and}\qquad 4/3\geq 1/2(1-\alpha)$$
that is
$$\alpha \geq 5/8\qquad \textrm{and} \qquad  \alpha\leq 5/8.$$

We therefore conclude that any $d_R$-geodesic between $X$ and $Y$ must
cross the $1$-simplex at the point $T_{5/8}$. The completely symmetric
calculation shows that any $d_L$-geodesic must cross the central edge
at the point $T_{3/8}$. Thus no path from $X$ to $Y$ can be
simultaneously $d_R$- and $d_L$-geodesics. It follows that no
$d$-geodesic in $CV_2$ joins $X$ and $Y$.

\section{Quasi-geodesics} \label{qg}
In section~\ref{s5} we have seen how to construct folding paths that
are $d_R$-geodesic. In this section we address the question of
whether such paths are quasi-geodesic for the symmetric metric, with
constants depending only on the rank. In other words, we ask whether
two points of outer space can be joined by a quasi-geodesic with
uniform constants.

To start, we recall the definition of a quasi-geodesic path.

\begin{dfn}
A path, $\alpha : I \to X$, where $I$ is a real interval and $(X, d)$
is a metric space, is called a $(\lambda, \epsilon)$ quasi-geodesic
if for every $x, y \in I$,
$$
 \frac{1}{\lambda} | x-y | - \epsilon \leq d(\alpha(x), \alpha(y))
 \leq \lambda |x-y | + \epsilon.
 $$
\end{dfn}

The following lemma is tautological.
\begin{lem}\label{l_3s3_7.1}
  Let $\alpha$ be a path from an interval to a metric space. Suppose
  that there is a constant $C$ such that
$$d(\alpha(x),\alpha(y))>C\cdot\operatorname{length}(\alpha|_{[x,y]}).$$
Then the arc-length reparameterisation of $\alpha$ is bi-lipschitzian
with constants $C,1$. In particular, it is a $(C,0)$ quasi-geodesic.
\end{lem}

\begin{thm}[4 point property]\label{t4pp}
 Let $A,B$ be two marked metric graph of the same rank. Let $\alpha$
 be a $d_R$-geodesic from $A$ to $B$ constructed as in
 Theorem~\ref{t3_drgeod}. Then for every $s\leq x\leq y\leq t$ we have
$$d(\alpha(s),\alpha(t))\geq d(\alpha(x),\alpha(y)).$$
\end{thm}
\proof Let us denote by $l_p$ the length function of the point
$\alpha(p)$. We consider the folding paths constructed before the
rescaling to volume 1, so that while volume is not constant along the
path, for every $p<q$, $\Lambda_R(\alpha(p), \alpha(q))=1$. Thus, the
distance between $\alpha(p)$ and $\alpha(q)$ is exactly the logarithm
of $\Lambda_L(\alpha(p), \alpha(q))=\sup l_p/l_q$. Now we look at the
points $s\leq x\leq y\leq t$. As in Proposition~\ref{p_optimal},
there exists a $\mu$ which realises $\Lambda_L(\alpha(x),
\alpha(y))$. Next we realise $\mu$ as an immersed path in
$\alpha(s)$. The folding path itself has two parts, one in which we
shrink the lengths of certain edges, and another in which we
isometrically identify edges - folding. In either of these parts it
is clear that the length of $\mu$ can never increase as we travel
along the path. Thus,
$$
l_s(\mu) \geq l_x(\mu) \geq l_y(\mu) \geq l_t(\mu).
$$
In particular,
$$
\sup l_s/l_t \geq l_s(\mu)/l_t(\mu) \geq l_x(\mu)/ l_y(\mu) = \sup
l_x/l_y, $$ and thus $d(\alpha(s),\alpha(t))\geq
d(\alpha(x),\alpha(y))$, as required. \qed

\begin{prop}\label{p4ppqg}
  Let $\gamma$ be a path with the 4 point property. Suppose that $\gamma$ is a
  finite union of pieces which are quasi-geodesic.
  Then $\gamma$ is a quasi-geodesic
  with constants depending on the constants of the pieces and on the
  number of the pieces.

  More precisely, if $\gamma$ is the path with the 4 point property which is the concatenation of
  $n$ $(\lambda, \epsilon)$ quasi-geodesics, then $\gamma$ is a $(n\lambda, n \epsilon)$ quasi-geodesic.
\end{prop}
\proof By hypothesis, there exist numbers $x_0 \leq x_1 \leq \ldots
\leq x_n$ such that $\gamma$ is a map from the interval $[x_0, x_n]$
and that each restriction, $\gamma|_{[x_i, x_{i+1}]}$ is a $(\lambda,
\epsilon)$ quasi-geodesic (we assume that $n > 1$ since otherwise
there is nothing to prove). Now consider $p \leq q \in [x_0, x_n]$,
and find $i, j$ such that $p \leq x_i \leq x_j \leq q$ so that $i$ is
minimal and $j$ is maximal (note that $i \geq 1$ and $j \leq n-1$).
It is clear that,
$$
\begin{array}{rcl}
d(\gamma(p), \gamma(q)) &  \leq &  d(\gamma(p), x_i) + \sum_{k=0}^{k=j-i-1}
d(x_{i+k}, x_{i+k+1}) + d(\gamma(x_j), q) \\
& \leq & \lambda (x_i-p) + \lambda \sum_{k=0}^{k=j-i-1} (x_{i+k+1}-x_{i+k})
\\&  & + \ \lambda (q-x_j) + (2+j-i)\epsilon \\
& \leq & \lambda (q-p) + n \epsilon.\\
\end{array}
$$
For the other inequality we note that, using the $x_r$, we have
divided the interval $[p,q]$ into at most $n$ pieces. Thus, one of
these pieces is of length at least $(q-p)/n$. Now, suppose that
$x_{i+k+1} - x_{i+k} \geq (q-p)/n$. Then, by the 4 point property,
$$
\begin{array}{rcl}
d(\gamma(p), \gamma(q)) &  \geq &  d(x_{i+k}, x_{i+k+1}) \\
&  \geq  &  (x_{i+k+1} - x_{i+k})/\lambda - \epsilon \\
& \geq &  (q-p)/n\lambda - \epsilon.\\
\end{array}
$$
Clearly, the same argument works if either $x_i - p \geq (q-p)/n$ or
$q - x_j \geq (q-p)/n$. \qed

\begin{example}
There are metric spaces with no rectifiable, non-constant paths
having the 4 point property.
\end{example}
\proof Consider the space $L^2([0,1])$ of the square-summable
functions on $[0,1]$. Let $f:[0,1]\to L^2([0,1])$ be the embedding
$$t\mapsto \chi_{[0,t]},$$
where $\chi_{[0,t]}$ denotes the characteristic function of the set
$[0,t]$. Let $d$ the $f$-pull-back metric on $[0,1]$:
$$d(s,t)=\sqrt{t-s}.$$
It is straightforward to check that $([0,1],d)$ has the 4 point
property and no rectifiable, non-constant paths.\qed

\medskip

By Theorem~\ref{t4pp} and Proposition~\ref{p4ppqg}, to check whether
a right geodesic between two points $A$ and $B$, constructed as in
Theorem~\ref{t3_drgeod}, is a quasi-geodesic (with uniform constants
not depending on $A$ and $B$,) it is enough to check whether the fast
folding path from $A_0$ to $\bar B$ is a quasi-geodesic.

\begin{dfn}[Multiplicities]
Let $A_t\neq B$ be any point in a fast folding path. The {\em
multiplicity} of a turn $\tau$ in a loop $\gamma$ is the number
$\mu_{\tau,t}(\gamma)$ of occurrences of $\tau$ turn in $\gamma$
(counted without any orientation.)

The {\em folding multiplicity} of $\gamma$ is the sum $\mu_t(\gamma)$
of the multiplicities of all folding turns (see
Definition~\ref{fastfold}) in $\gamma$:
$$\mu_t(\gamma)=\sum_\tau \mu_{\tau,t}(\gamma).$$
 \end{dfn}

In order to use Lemma~\ref{l_3s3_7.1}, we need to estimate the local
speed of a fast folding path. A folding path is PL, and therefore
smooth in all but finitely many points (w.r.t. the PL-structure of
$CV_n$.) In particular, the right-derivative is always defined, and
its integral gives the total length of the path.

\begin{lem}[Local speed of a folding path]\label{l3s3_7.6}
  Let $t\mapsto A_t$ be a fast folding path. Then, its local speed is
$$\frac{2\mu_t(\gamma)}{l_{A_t}(\gamma)}$$
where $\gamma$ is a folded loop minimising
$l_{A_t}(\gamma)/\mu_t(\gamma)$.
\end{lem}
\proof Recall that in our situation (isometric folding as in
Theorems~\ref{t3_drgeod},) we have $d=d_L$. Therefore, for small
enough $\varepsilon$, the distance between $A_{t+\varepsilon}$ and
$A_t$ is given by
$$d(A_{t+\varepsilon},A_t)
=\log(\sup_{\xi}\frac{l_{A_t}(\xi)}{l_{A_{t+\varepsilon}}(\xi)})
=\log(\sup_{\xi}\frac{l_{A_t}(\xi)}{l_{A_t}(\xi)-2\mu_t(\xi)\varepsilon})
$$
which is thus realised by a loop $\gamma$ minimising
$l_{A_t}(\gamma)/\mu_t(\gamma)$. Note that $\gamma$ can be always
chosen to be simple.

Therefore, the speed (as right-derivative) is given by
$$\lim_{\varepsilon \to 0}\frac{d(A_{t+\varepsilon},A_t)}{\varepsilon}
=\lim_{\varepsilon \to 0}
\frac{1}{\varepsilon}\log(\frac{l_{A_t}(\gamma)}{l_{A_{t+\varepsilon}}(\gamma)})
=\lim_{\varepsilon \to 0} \frac{1}{\varepsilon}
\log(\frac{l_{A_t}(\gamma)}{l_{A_t}(\gamma)-2\mu_t(\gamma)\varepsilon})
=\frac{2\mu_t(\gamma)}{l_{A_t}(\gamma)}.
$$
\qed

Another quantity we need to estimate during a folding procedure, is
the speed we are approaching the final point $B$, defined as the
right-derivative of the distance from $B$.

\begin{lem}[Local speed toward $B$]\label{l3s3_7.7}
 Let $t\mapsto A_t$ be a fast folding path. Then, the speed at which $A_t$
 is approaching $B$ is given by
$$\frac{2 \mu_t(\gamma)}{l_{A_t}(\gamma)}$$
where $\gamma$ is a loop that realises the maximal stretching factor
from $B$ to $A_t$.
\end{lem}
\proof As above, since $t\mapsto A_t$ is an isometric folding path
constructed as in Theorem~\ref{t3_drgeod}, we are interested only in
$d_L$. We have
$$d(A_t,B)=d_L(A_t,B)=\log(\frac{l_{A_t}(\gamma)}{l_B(\gamma)}).
$$
During  the folding procedure, in the marked graph $A_t$, the length
of $\gamma$ decrease twice the number of occurrences of the folding
turns in $\gamma$. Whence the claim follows.\qed

\medskip

Now, the aim is to show that the ratio between the speed toward $B$
and the local speed is bounded below by a given constant. Indeed, if
so, one could deduce that the hypothesis of Lemma~\ref{l_3s3_7.1} is
satisfied, this providing quasi-geodesics with uniform constants.

\begin{lem}
Let $A_t\neq B$ be any point in a fast folding path.
 Let $\gamma$ be a loop that
realises the maximal stretching factor from $B$ to $A_t$. Then
$$\mu_t(\gamma)\geq 1$$
\end{lem}

\proof Otherwise $\gamma$ would be  immersed via the optimal map $f$
used for defining the folding procedure, which would imply
$l_{A_t}(\gamma)=l_B(\gamma)$, whence $A_t=B$.\qed

\begin{lem}\label{l3s5_7.13}
In a fast folding path, for any loop $\gamma$, the quantity
$\mu_t(\gamma)$, as a function
  of $t$, is monotone non-increasing.
\end{lem}
\proof Let $0=t_0<t_1...$ be the subdivision of times. Clearly,
nothing change for $t$ different from the $t_i$'s. We show that the
multiplicity cannot increase passing trough any $t_i$'s. Let
$\tau=(a,b)$ be a turn where the segments $a$ and $b$ are identified
during the interval of time $[t_{i-1},t_i]$. The segments $a$ and $b$
have one extreme in common, say the starting point. On the other
hand, the ending points of $a$ and $b$, say $x$ and $y$ respectively,
must be different, otherwise the folding procedure would decrease the
rank of our marked metric graphs, which is not possible.

The multiplicity, in $\gamma$, of the turns that already exist for
$t\in(t_{i-1},t_i)$ is unchanged. So we have to check what happens to
the new turns created by the folding. Those are pair of segments $a'$
and $b'$ having $x$ and $y$ as starting points, and identified by the
optimal map. Let $\{(a_j,b_j)\}$ be the set of turns folded for
$t\in(t_{i-1},t_i)$ whose ending points are $x$ and $y$.

The multiplicity of the turn $(a',b')$ counts how many times $\gamma$
passes trough the turn. But any times that $\gamma$ passes trough
$(a',b')$ must passes trough one of the $(a_j,b_j)$'s as well. So the
total sum is not increased.\qed

\begin{dfn}[$\varepsilon$-thin part]
  The {\em $\varepsilon$-thin part} of $CV_n$ is the set of marked metric
  graphs having a loop shorter than $\varepsilon$ in the
  volume-one-representative. In other words, the class
 a marked metric graph $A$ in $CV_n$ lies in the $\varepsilon$-thin
 part if
$$\frac{l_A(\textrm{shortest loop of }A)}{\operatorname{vol}
  A}<\varepsilon.$$
Otherwise, we say that $A$ lies in the {\em $\varepsilon$-thick
part}.

\end{dfn}

\begin{lem}\label{l3s5_7.15}
There is a constant $C>0$ such that for any fast folding path
$t\mapsto A_t$, if $A_t$ never enters the $\varepsilon$-thin part,
then the ratio between the speed approaching toward $B$ and the local
speed is bounded below by $C\cdot\varepsilon$.
\end{lem}
\proof Since our folding procedure is isometric, if, starting from
$A_t$, we fold during a time $T$, then the volume of $A_t$ is
decreased at least by $T$:
$$T\leq\operatorname{vol}(A_t)-\operatorname{vol}(B)=
\operatorname{vol}(A_t)-1.$$
 On the other
hand, the length of a given loop is decreased by
$$l_{A_t}(\gamma)-l_B(\gamma)=2\int_t^{t+T}\mu_s(\gamma)\,ds
\leq 2T\mu_t(\gamma)
$$
where the inequality follow from Lemma~\ref{l3s5_7.13}.

Now, let $\gamma$ be a loop realising the maximal stretching factor
from $B$ to $A_t$. Since $\operatorname{vol} (\bar B)=1$, the length
of $\gamma$ in $\bar B$ is less than $2$ (because of
Proposition~\ref{glasses}.) By the above inequalities it follows that

$$l_{A_t}(\gamma)\leq 2\mu_t(\gamma)\operatorname{vol}(A_t).$$

Let $\gamma_1$ be simple a loop minimising $l_{A_t}(w)/\mu_t(w)$. The
ratio between the approaching speed toward $B$ and the local speed
is, by Lemmata~\ref{l3s3_7.6} and~\ref{l3s3_7.7}
$$
\frac{l_{A_t}(\gamma_1)\mu_t(\gamma)}{l_{A_t}(\gamma)\mu_t(\gamma_1)}
$$
which is therefore bounded below by
$$
\frac{l_{A_t}(\gamma_1)}{2\operatorname{vol}{A_t}\mu_t(\gamma_1)}\geq
C\frac{l_{A_t}(\textrm{shortest loop of
}A_t)}{\operatorname{vol}(A_t)}
$$

where $C$ is a constant depending only on the rank $n$. Actually, the
constant $C$ depends on the fact that $\mu_t(\gamma_1)$ is bounded
above, depending on the rank, because $\gamma_1$ is a simple loop.

Therefore, the ratio between the approaching speed toward $B$ and the
local speed is bounded below by $C\cdot\varepsilon$ if $A_t$ lies in
the $\varepsilon$-thick part of $CV_n$.\qed

\medskip

An immediate corollary is the following
\begin{thm}[Folding paths are quasi-geodesic]\label{t_4s4_6.13}
  For any $\varepsilon>0$
  there are constants $K,L$ depending only on $\varepsilon$ and the
  rank of $CV_n$ such that for any two marked
  metric graph $A$ and $B$ whose corresponding fast folding path
  $t\mapsto A_t$ from $A_0$ to $\bar B$ (notation as in
  Theorem~\ref{t3_drgeod}) stay in the $\varepsilon$-thick part, there
  is right-geodesic between $A$ and $B$ which is a $(K,L)$-quasi-geodesic.
\end{thm}

\proof Lemma~\ref{l3s5_7.15} implies that the hypothesis of
Lemma~\ref{l_3s3_7.1} is satisfied. By Theorem~\ref{t4pp} and
Proposition~\ref{p4ppqg} the claim follows.\qed

\section{Iterating automorphisms}
\label{iterate}
Here, we study the behaviour of the orbits of automorphisms with
respect to our metrics.

\begin{thm}
\label{orbit}
  Let $\Phi\in\operatorname{Aut}(F_n)$ be an automorphism of
  exponential growth. Then for any $A\in CV_n$ the sequence $\Phi^hA$
  is a quasi-geodesic as a map from $\mathbb Z\to CV_n$. Moreover, if
  $A$ is a train-track for $\Phi$, then it is a $d_R$-geodesic.
\end{thm}

\proof

 If $\Phi$ has exponential growth so does $\Phi^{-1}$ (this is a
 consequence of the existence of the relative train track
 representatives of~\cite{MR1147956}.)  That means that $\sup_{1\neq w\in F_n}
l(\Phi^h(w))/l(w)>kc^h$ for some $k>0$ and $c>1$, where the length
$l$ is calculated in any fixed rose
 (and the same holds for $\Phi^{-1}$.)
We have
$$\sup_{1\neq w\in F_n}\frac{l_A(\Phi^{h+m}w)}{l_A(\Phi^mw)}=
\sup_{1\neq w\in F_n}\frac{l_A(\Phi^h w)}{l_A(w)} =\sup_{1\neq w\in
F_n}\frac{l_A(\Phi^w w)}{l(\Phi^h w)} \cdot\frac{l(\Phi^h)}{l(w)}
\cdot\frac{l(w)}{l_A(w)}
$$

In the last term of above inequality, the first and the last factors
are bounded below by constants because $A$ lies at finite distance
from the rose used for calculating $l$. The middle term is bounded
below by $kc^h$ by our hypothesis of exponential growth. Similarly,
using that also $\Phi^{-1}$ has exponential growth, we can show that
$$\Lambda(\Phi^{h+m}A,\Phi^mA)>kc^h$$
for some constants $k>0$ and $c>1$, this giving
$$d(\Phi^{h+m}A,\Phi^mA)>\log k+h \log c$$

The other inequality is even easier, and does not need any assumption
on $\Phi$:
$$\sup_{1\neq w\in F_n}\frac{l_A(\Phi^{h+m}w)}{l_A(\Phi^mw)}=
\sup_{1\neq w\in F_n}
\frac{l_A(\Phi^{l+m}w)}{l_A(\Phi^{h+m-1}w)}\cdot
\frac{l_A(\Phi^{h+m-1}w)}{l_A(\Phi^{h+m-2}w)}\cdots
\frac{l_A(\Phi^{1+m}w)}{l_A(\Phi^mw)}$$ which is bounded above by
$$\left(\sup_{1\neq w\in F_n}\frac{l_A(\Phi w)}{l_A( w)}\right)^h$$
whence (arguing the same way for $\Phi^{-1}$)
$$\Lambda(\Phi^{h+m}A,\Phi^m A)\leq\Lambda(\Phi A,A)^h$$
and
$$d(\Phi^{h+m}A,\Phi^m A)\leq hd(\Phi A,A).$$

Suppose now that $A$ is a train track for $\Phi$. Then every edge is
stretched exactly by $\lambda$, the Perron-Frobenius eigenvalue
associate to the transition matrix for $\Phi$ (see~\cite{MR1147956}.)
It follows that $\Lambda_R(\Phi^{h+m}, \Phi^m)=\lambda^h$, and the
second claim follows.\qed

\medskip

The fact that train tracks for $\Phi$ and $\Phi^{-1}$ are in general
different, and that also the Perron-Frobenius eigenvalues for $\Phi$
and $\Phi^{-1}$ may differ, tells us that we cannot follows this
approach for building a $d$-geodesic axis for $\Phi$.

Now, Theorem~\ref{orbit} clearly fails if the automorphism in
question is of polynomial growth. However it is important to note
that, nevertheless, the various folding paths from a point to the
points in its orbit may still be quasi-geodesics (with the unit speed
parametrisation) as in the following example.

\begin{example}
Let $R$ be the rose of rank 2, with loops labelled $A, B$ and let
$\phi$ be the automorphism which sends $A$ to $A$ and $B$ to $BA$.
Then, for any $k$, the folding path from $R$ to $\phi(R)$ is a
$(4,0)$ quasi-geodesic.
\end{example}
\proof In the rose, the petals have the same length, but since our
metric is scale invariant, we may choose that length - we choose it
to be $k+1$. We let $R_k$ denote $\phi^k(R)$, which then also has two
loops of the same length, which we label $A_k$ and $B_k$, and give
them both length $1$. By definition, $A$ maps to the loop $A_k$ in
$R_k$ and $B$ maps to $B_k (A_k)^k$.

In the folding path we start, first of all, by shrinking all the
edges so that (after scaling, which we have already done) the map
from the left to the right is isometric on edges. This means that we
shrink the loop $A$ until it has length $1$. We call this new graph
$R_0$; it has one vertex and two loops, $A_0 \to A_k$ and $B_0 \to
B_k (A_k)^k$. The length of $A_0$ is $1$ and the length of $B_0$ is
$k+1$.

The folding path then proceeds by folding $A_0$ into $B_0$. If one
imagines this as a discrete process, after the $i^{th}$ stage we will
obtain a graph $R_i$, with a single vertex and two loops, $A_i \to
A_k$ and $B_i \to B_k (A_k)^{(k-i)}$; the length of $A_i$ is $1$ and
the length of $B_i$ is $K+1-i$.

If we then fold a part of $A_i$, of length $\delta$, into $B_i$ we
travel to a point in the folding path which we shall call $R_{i,
\delta}$. This has two vertices, $\bullet$ and $\circ$, and three
edges, $A_{i, \delta}, B_{i, \delta}$ and $C_{i, \delta}$.

\begin{figure}[htbp]
\begin{center}
\ \includegraphics[width=3in]{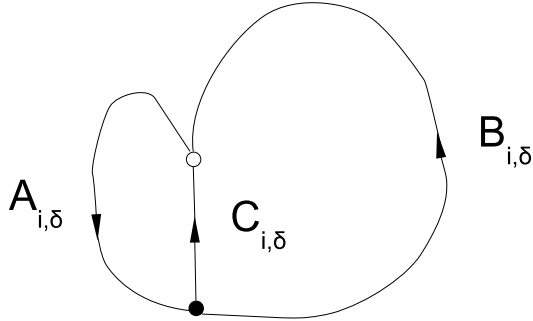}
\end{center}
\caption{The graph $R_{i, \delta}$} \label{graph}
\end{figure}

Here, we can map the vertex $\bullet$ to the unique vertex of $R_k$
and then map the loop $A_{i, \delta} C_{i, \delta}$ to $A$ and $B_{i,
\delta} C_{i, \delta}$ to $B_k (A_k)^{(k-i)}$ (this is enough to
specify the marking up to homotopy equivalence); the length of $A_{i,
\delta}$ is $1-\delta$, the length of $B_{i, \delta}$ is
$k+1-i-\delta$ and the length of $C_{i,\delta}$ is $\delta$. This
marked metric graph represents an arbitrary point on the folding path
from $R_0$ to $R_k$. Now, following Lemma~\ref{l3s3_7.6}, the local
speed is realised by the loop $B_{i, \delta} \overline{A_{i,
\delta}}$, whereas the distance to $R_k$ is realised by the loop
$B_k$ which is realised by $B_{i, \delta} \overline{A_{i, \delta}} (
\overline{C_{i, \delta}} \ \overline{A_{i, \delta}})^{k-i-1}$ in
$R_{i, \delta}$. Both of these loops pass through the unique folding
turn (Definition~\ref{fastfold}) of $R_{i, \delta}$ exactly once.

Hence, by Lemmas~\ref{l3s3_7.6} and \ref{l3s3_7.7}, the ratio of the
speed toward $R_k$  and the local speed is,
$$
\frac{k+2-i-2\delta}{2k+1-2i-2\delta} \geq \frac{1}{2}.
$$
Thus, by Lemma~\ref{l_3s3_7.1}, the path from $R_0$ to $R_k$ is a
$(2,0)$ quasi-geodesic and thus by Proposition~\ref{p4ppqg}, the
whole path is a $(4,0)$ quasi-geodesic. \qed

\section{Some open questions}

In this section we address some questions which arose during the many
conversations we had with colleagues, principally during
the coffee breaks of conferences, about the metric properties
of Outer Space.

\subsection{Existence of quasi-geodesics.} As we've seen, folding paths
that do not fold into the thin part provide quasi-geodesics for the
symmetric metric. Here we address mainly two questions. First, whether
a folding path will always produce a quasi-geodesic or not, with constants depending only on the rank. Second, whether it is in general possible
to connect any two marked metric graphs with a path which is a
quasi-geodesic, with constants depending only on the rank of the
graphs.

For the latter question, there is an heuristic argument:
suppose the answer is no. Then, letting blowing up the constants, one
would get a counter-example-sequence that contradicts
Lemma~\ref{l_3s3_7.1}. Then, following the arguments of
Theorem~\ref{t_4s4_6.13} one gets that the folding paths of the
counter-example-sequence  will
eventually enter any $\varepsilon$-thin part, but explicit
computations show that a folding path that enters the thin part cannot
stay for too long inside that part (one has perhaps to understand how many times a folding path can enter the thin part.) Thus suggesting an
affirmative answer to our questions.

\subsection{Existence of a geodesic axis for an iwip.}
We have seen that iterates of automorphisms produce quasi-geodesics (and
geodesics for the non-symmetric metric.) The natural question here is
whether automorphisms have an axis and whether can such an axis be
described in terms of metric properties. Also, one can ask whether one
can compute the ``geometric rank'' of such axis. Is there any
analogue of the bounded projection Lemma? (see~\cite{BF1} and the
recent preprint~\cite{Ya}.)

\subsection{Hyperbolicity, flats and coarse properties}
It is natural to ask whether some subset of Outer space (some
thick-part?) is hyperbolic or presents hyperbolicity phenomena. On the
other hand, it would be interesting to study the (quasi-) flats of
Outer space, if any. In general coarse properties of Outer Space are
still unknown (for instance, what do its asymptotic cones look like?)

\end{document}